\documentclass[10pt]{article}
\usepackage{amssymb}
\usepackage{graphicx}
\usepackage{amsmath}

\newtheorem{theorem}{Theorem}

\newtheorem{problem}[theorem]{Problem}
\newtheorem{proposition}[theorem]{Proposition}
\newtheorem{remark}[theorem]{Remark}

\newenvironment{proof}[1][Proof]{\textbf{#1.} }{\ \rule{0.5em}{0.5em}}
\input{tcilatex}

\begin{document}

\title{ Computational Topology of Equivariant Maps\\
 from Spheres to Complements of Arrangements}
\author{Pavle V. M. Blagojevi\'{c}, Sini\v{s}a\ T. Vre\'{c}ica, Rade T. \v{Z}%
ivaljevi\'{c}}
\date{March 2004}
\maketitle

\section{{Introduction}}

\subsection{Computational Equivariant Topology}

\label{sec:CET}

\noindent Perhaps one of the main general (unsolved) problems of
Computational Topology is to determine the scope of the field and major
lines of prospective research. This includes the identification of classes
of \textquotedblleft model problems\textquotedblright\ where topology and
computational mathematics interact in an essential way. The following
references may serve as a source of initial information about Computational
Topology and as a guide for some of the existing applications.

\begin{enumerate}
\item[(ATCS)] G.~Carlsson (ed.) \textit{Proceedings of the conference on
Algebraic Topological Methods in Computer Science}. To appear in the journal
Homology, Homotopy and Applications.

\item[(CTOP)] M.~Bernet. al. Emerging challenges in computational topology.
\textit{ACM Computing Research Repository}. {arXive: cs.CG/9909001}.

\item[(TM)] R.~{\v{Z}}ivaljevi\'{c}. Topological Methods. Chapter 14 of
\textit{Handbook of Discrete and Computational Geometry} (J.E.~Goodman,
J.~O'Rourke, eds.), new edition, CRC Press, Boca Raton 2004.
\end{enumerate}

\noindent One of our objectives in this paper is to identify the problem of
calculating the \textit{topological obstructions for the existence of
equivariant maps} as one of the problems paradigmatic for computational
topology. Recall that, given a group $G$, a $G$-equivariant map $%
f:X\rightarrow Y$ between two $G$-spaces is a symmetry preserving map, i.e.
a map satisfying the condition $f(g\cdot x)=g\cdot f(x)$. In applications in
discrete and computational geometry, $X$ is usually a manifold of all
\textquotedblleft feasible\textquotedblright\ configurations (the
configurations space) while $Y$ is typically a complement of a real, affine
subspace arrangement (the test space), see (TM).

\noindent The problem of calculating/evaluating the complexity of
equivariant obstructions has many aspects relevant for computational
topology, including the following.

\begin{enumerate}
\item[($A_{1}$)] The (non)existence of an equivariant map is an essential
ingredient in the application of the \textit{configuration space/test map}%
-scheme (see (TM)), which has proven to be a very effective tool in solving
combinatorial or discrete geometric problems of relevance to computing and
analysis of algorithms, \cite{Mat94} \cite{guide2}.

\item[($A_{2}$)] The existence of an equivariant map can often be
interpreted as a problem of mapping an object (often a cell of some
dimension) to an Euclidean space, subject to some boundary constraints and
avoiding some obstacles (arrangement of subspaces). This aspect can be seen
as a relative of the motion planning problem from robotics (the case of a $1$%
-dimensional cell). Recall that the existence of an equivariant map is
equivalent to the problem of sectioning a (vector) bundle and the condition
of avoiding obstacles is translated into the question of sectioning a bundle
subject to some additional constraints (relations). The latter problem is of
utmost importance in many areas of mathematics, including Gromov convex
integration theory, combinatorial geometry on vector bundles \cite{T-V} etc.
\end{enumerate}

\noindent

\noindent \textbf{Conclusion:} Both positive and negative aspect of the
problem of the existence of an equivariant map is of theoretical interest.
In the negative case, i.e. if an equivariant map exists, one should be able
to describe and evaluate the complexity of an algorithm that constructs such
a map. One way to achieve this goal is to develop an \textquotedblleft
effective\textquotedblright\ obstruction theory. This is of particular
interest in the context of ($A_{2}$) where such a theory would allow
effective placements of objects into an environment with obstacles, subject
to some boundary constraints.

\bigskip

\noindent \noindent \textbf{This paper: }We use a well known problem in
discrete and computational geometry (Problem 1) as a motivation and an
associated question \textbf{\ }from equivariant topology (Problem 3) as a
point of departure to illustrate many aspects, both theoretical and
computational, of the general \emph{existence of en equivariant map }%
problem. A variety of techniques are introduced and discussed with the
emphasis on concrete and explicit calculations. This eventually leads
(Theorems 18 and 19) to an almost exhaustive analysis of when such maps do
or do not exist in this particular case of interest.\bigskip

\subsection{Obstructions to the existence of equivariant maps}

\noindent Suppose that $X$ and $Y$ are topological spaces and let $G$ be a
groups acting on both of them. Suppose that $Z$ is a closed, $G$-invariant
subspace of $Y$. We focus our attention in this paper on spaces which are
simplicial or $CW$-complexes and finite groups $G$, see \cite{Mat94} or (TM)
for a glossary of basic topological terms. \noindent The problem of deciding
if there exists a $G$-equivariant map $f:X\rightarrow Y\setminus Z$ can be
approached by two, closely related ways. The first is to build such a map
step by step, defining it on the skeletons $X^{(k)}$ one at a time, and
attempting to extend it to next skeleton $X^{(k+1)}$. The second is to start
with a sufficiently generic, equivariant map $f:X\rightarrow Y$ and, in the
case the singularity $S(f):=f^{-1}(Z)$ is nonempty, try to modify $f$ in
attempt to make the set $S(f)$ vanish. In both approaches there may appear
\emph{obstructions} which prevent us from completing the process, thus
showing that such equivariant maps $f:X\rightarrow Y\setminus Z$ do not
exist. In the first approach the obstruction is evaluated in the
corresponding \emph{equivariant cohomology group}, while in the second it
lies in a \textquotedblleft dual\textquotedblright\ \emph{equivariant
homology group}, see Section~\ref{sec:obs-theo} for some technical details
and references.

\section{{The motivating problem}}

\noindent Many aspects of \textit{the existence of equivariant map problem}
relevant for computational topology are particularly well illustrated by the
problem of finding $k$-fan partitions of spherical measures, see Problem~\ref%
{OsnovniProblem} and original references \cite{BaMa2001}, \cite{BaMa2002}
and \cite{VreZiv2002}. Recall that this question arose in connection with
some partition problems in discrete and computational geometry, \cite%
{Aki2000}, \cite{Bes2000}, \cite{Ito2000}, \cite{Kan1999}, \cite{Sak1998}
\cite{VreZiv2001}. It turns out that Problem~\ref{OsnovniProblem} is closely
related to a problem of the existence of equivariant maps, see Problem~\ref%
{prob:dihedral} in Section~\ref{sec:scheme}.

\noindent In this paper we give a fairly complete analysis, Theorems~{\ref%
{Th.Main1}} and {\ref{Th.Main2}}, of the Problem~\ref{prob:dihedral} and
solve it in all but a few exceptional cases.

\subsection{Partition of measures by $k$-fans}

\noindent A $k$-fan $\mathfrak{p}=(x;l_{1},l_{2},\ldots ,l_{k})$ on the
sphere $S^{2}$ is a point $x$, called the center of the fan, and $k$ great
semicircles $l_{1},\ldots ,l_{k}$ emanating from $x$. We always assume
counter clockwise enumeration of great semicircles $l_{1},\ldots ,l_{k}$ of
a $k$-fan. Sometimes we use the notation $\mathfrak{p}=(x;\sigma _{1},\sigma
_{2},\ldots ,\sigma _{k})$, where $\sigma _{i}$ denotes the open angular
sector between $l_{i}$ and $l_{i+1},\,i=1,\ldots ,k$.

\noindent Let $\mu _{1},\mu _{2},\ldots ,\mu _{m}$ be \textit{proper} Borel
probability measures on $S^{2}$. Measure $\mu $ is \textit{proper} if $\mu
([a,b])=0$ for any circular arc $[a,b]\subset S^{2}$ and $\mu (U)>0$ for
each nonempty open set $U\subset S^{2}$. All the results can be extended to
more general measures, including the counting measures of finite sets, see
\cite{BaMa2001}, \cite{elisat}, \cite{T-V} for related examples.

\noindent Let $(\alpha _{1},\alpha _{2},\ldots ,\alpha _{k})\in \mathbb{R}%
_{>0}^{k}$ be a vector where $\alpha _{1}+\alpha _{2}+\ldots +\alpha _{k}=1$%
. Following \cite{BaMa2001}, and keeping in mind that we deal only with
proper measures, we say that a $k$-fan $(x;l_{1},\ldots ,l_{k})$ is an $%
\alpha $\textit{-partition} for the collection $\{\mu _{j}\}_{j=1}^{m}$ of
measures if
\begin{equation*}
(\forall i=1,\ldots ,k)\,(\forall j=1,\ldots ,m)\,\mu _{j}(\sigma
_{i})=\alpha _{i}.
\end{equation*}%
As in \cite{VreZiv2002}, a vector $\alpha \in \mathbb{R}^{k}$ is called $%
(m,k)$-\textit{admissible}, if for any collection of $m$ (proper) measures
on $S^{2}$, there exists a simultaneous $\alpha $-partition. The collection
of all $(m,k)$-admissible vectors is denoted by $\mathcal{A}_{m,k}$. Here is
the central problem about partitions of measures by $k$-fans.

\begin{problem}[\protect\cite{BaMa2001}, \protect\cite{BaMa2002},
\protect\cite{VreZiv2002}]
\ \label{OsnovniProblem}Describe the set $\mathcal{A}_{m,k}$ or
equivalently, find integers $m$, $k$ and vectors $\alpha \in \mathbb{R}^{k}$
such that for any collection $\mathcal{M}=\{\mu _{1},\mu _{2},\ldots ,\mu
_{m}\}$ of $m$ (proper) measures , there exist an $\alpha $-partition for $%
\mathcal{M}$.
\end{problem}

\noindent The analysis given in \cite{BaMa2001} shows that the most
interesting cases are $(3,2),(2,3),(2,4)$. It was shown (\cite{BaMa2001},
\cite{VreZiv2002}) that $\{(\frac{1}{2},\frac{1}{2}),(\frac{1}{3},\frac{2}{3}%
),(\frac{2}{3},\frac{1}{3})\}\subseteq \mathcal{A}_{3,2}$ and $\{(\frac{1}{4}%
,\frac{1}{4},\frac{1}{4},\frac{1}{4}),\,(\frac{1}{5},\frac{1}{5},\frac{1}{5},%
\frac{2}{5})\}\subseteq \mathcal{A}_{2,4}$.

\subsection{The configuration space / test map scheme}

\label{sec:scheme} The configuration space/test map scheme is a fairly
general method of translating combinatorial geometric problems into
topological problems, more precisely problems involving equivariant maps,
see (TM) or references \cite{Mat94}, \cite{elisat}, \cite{guide2}. It was
demonstrated in \cite{BaMa2001} that the problem of $\alpha $-partitions of
spherical measures also admits such a translation. Recall the key steps of
this elegant construction.

\noindent Let $\mu $ and $\nu $ be two Borel, probability measures on $S^{2}$
and $F_{k}$ the space of all $k$-fans on the sphere $S^{2}$. The space $%
X_{\mu }$ of all $n$-equipartitions of the measure $\mu $ is defined by
\begin{equation*}
X_{\mu }=\{(x;l_{1},\ldots ,l_{n})\in F_{n}\mid (\forall i=1,\ldots ,n)\,\mu
(\sigma _{i})=\tfrac{{1}}{{n}}\}.
\end{equation*}%
Observe that every $n$-fan $(x;l_{1},\ldots ,l_{n})\in X_{\mu }$ is
completely determined by the pair $(x,l_{1})$ or equivalently the pair $%
(x,y) $, where $y$ is the unit tangent vector to $l_{1}$ at $x$. Thus, the
space $X_{\mu }$ is Stiefel manifold $V_{2}(\mathbb{R}^{3})$ of all
orthonormal $2$-frames in $\mathbb{R}^{3}$. Recall that $V_{2}(\mathbb{R}%
^{3})\cong SO(3)\cong \mathbb{R}P^{3}$.

\noindent Let $\mathbb{R}^{n}$ be an Euclidean space with the standard
orthonormal basis $e_{1},e_{2},\ldots ,e_{n}$ and the associated coordinate
functions $x_{1},x_{2},\ldots ,x_{n}$. Let $W_{n}$ be the hyperplane $\{x\in
\mathbb{R}^{n}\mid x_{1}+x_{2}+\ldots +x_{n}=0\}$ in $\mathbb{R}^{n}$ and
suppose that the $\alpha $-vector has the form $\alpha =(\frac{a_{1}}{n},%
\frac{a_{2}}{n},\frac{a_{3}}{n},\frac{a_{4}}{n})\in \frac{1}{n}\,\mathbb{N}%
^{4}\subset \mathbb{Q}^{4}$ where $a_{1}+a_{2}+a_{3}+a_{4}=n$. Then the test
map $F_{\nu }:X_{\mu }\rightarrow W_{n}\subset \mathbb{R}^{n}$ for the
measure $\nu $ is defined by
\begin{equation*}
F_{\nu }(x,y)=F_{\nu }(x;l_{1},\ldots ,l_{n})=(\nu (\sigma _{1})-\tfrac{{1}}{%
{n}},\nu (\sigma _{2})-\tfrac{{1}}{{n}},\ldots ,\nu (\sigma _{n})-\tfrac{{1}%
}{{n}}).
\end{equation*}

\noindent The natural group of symmetries arising in this problem is the
dihedral group $D_{2n}$. It is interesting and sometimes useful to forget a
part of the structure and study the associated problem relative to other
subgroups $G\subset D_{2n}$. The associated test subspace is in that case
the space $D_{G}(\alpha )=D(\alpha ):=\cup \mathcal{A}(\alpha )\subset W_{n}$
defined as the union of the smallest $G$-invariant linear subspace
arrangement $\mathcal{A}(\alpha )$ in $\mathbb{R}^{n}$ containing the linear
subspace $L(\alpha )\subset W_{n}$ defined by
\begin{equation}
L(\alpha ):=\{x\in \mathbb{R}^{n}\mid
z_{1}(x)=z_{2}(x)=z_{3}(x)=z_{4}(x)=0\},  \label{L(alfa)}
\end{equation}%
where
\begin{equation*}
\begin{array}{lcl}
z_{1}(x)=x_{1}+x_{2}+\ldots +x_{a_{1}}, &  & z_{2}(x)=x_{a_{1}+1}+\ldots
+x_{a_{1}+a_{2}}, \\
z_{3}(x)=x_{a_{1}+a_{2}+1}+\ldots +x_{a_{1}+a_{2}+a_{3}}, &  &
z_{4}(x)=x_{a_{1}+a_{2}+a_{3}+1}+\ldots +x_{n}.%
\end{array}%
\end{equation*}%
Of course, our central interest is in the case $G=D_{2n}$. Recall \cite%
{BaMa2001} that the action of $D_{2n}=\langle \omega ,\varepsilon
\,|\,\omega ^{n}=\varepsilon ^{2}=1,\,\omega \varepsilon =\varepsilon \omega
^{n-1}\,\rangle $ on the configuration space $X_{\mu }$ and the test space $%
W_{n}$ respectively, is given by
\begin{equation*}
\begin{array}{lcl}
\omega (x;l_{1},\ldots ,l_{n})=(x;l_{2},\ldots ,l_{n},l_{1}), &  &
\varepsilon (x;l_{1},\ldots ,l_{n})=(-x;l_{1},l_{n},l_{n-1,}\ldots ,l_{2}),
\\
\omega (x_{1},\ldots ,x_{n})=(x_{2},\ldots ,x_{n},x_{1}), &  & \varepsilon
(x_{1},\ldots ,x_{n})=(x_{n},\ldots ,x_{2},x_{1}).%
\end{array}%
\end{equation*}%
where $(x;l_{1},\ldots ,l_{n})\in X_{\mu }$ and $(x_{1},\ldots ,x_{n})\in
W_{n}$. $D_{2n}$ acts also on the complement $M(\alpha )=W_{n}\setminus
D(\alpha )$ since $\mathcal{A}(\alpha )$ is the smallest $D_{2n}$-invariant,
linear subspace arrangement which contains linear subspace $L(\alpha )$.
Note that $D_{2n}$-action on the configuration space $X_{\mu }$ is free.
Obviously the test map $F_{\nu }$ is a $D_{2n}$-equivariant map. Note that
the configuration space $X_{\mu }$ is $D_{2n}$-homeomorphic to the Stiefel
manifold $V_{2}(\mathbb{R}^{3})$, which as a $D_{2n}$-space has the action
described by
\begin{equation*}
\omega (x,y)=(x,R_{x}({\frac{2\pi }{n}})(y))\text{, }\varepsilon (x,y)=(-x,y)
\end{equation*}%
where $R_{x}(\theta ):\mathbb{R}^{3}\rightarrow \mathbb{R}^{3}$ is the
rotation around the axes determined by $x$ through the angle $\theta $. The
following proposition is a consequence of this analysis.

\begin{proposition}
\label{prop:albeit} Let $\alpha =(\frac{a_{1}}{n},\frac{a_{2}}{n},\frac{a_{3}%
}{n},\frac{a_{4}}{n})\in \frac{1}{n}\,\mathbb{N}^{4}\subset \mathbb{Q}^{4}$
be a vector such that $a_{1}+a_{2}+a_{3}+a_{4}=n$. Let $G$ be a subgroup of
the dihedral group $D_{2n}$. If there does not exist a $G$-map $F:V_{2}(%
\mathbb{R}^{3})\rightarrow M_{G}(\alpha )$, where $M_{G}(\alpha
):=W_{n}\setminus D_{G}(\alpha )$, then for any two measures $\mu $ and $\nu
$ on $S^{2}$, there always exists a $4$-fan which simultaneously $\alpha $%
-partitions both $\mu $ and $\nu $. In other words the non-existence of such
a map implies $\alpha \in \mathcal{A}_{2,4}$.
\end{proposition}

\noindent Recall that in this paper the emphasis is put on computational
equivariant topology and the question of the existence of equivariant maps.
Hence, in light of Proposition\ref{prop:albeit}, \ it is natural to focus
our attention on the following problem as a problem closely related, albeit
not equivalent, to the initial $\alpha $-partition problem.

\begin{problem}
\label{prob:dihedral} For a given subgroup $G$ of the dihedral group $D_{2n}$%
, and the associated $G$-spaces $V_{2}(\mathbb{R}^{3})$ and $M_{G}(\alpha )$%
, find an explicit description of the set
\begin{equation*}
\Lambda _{G}=\{\alpha \in \tfrac{1}{n}\,\mathbb{N}^{4}\,|\,\text{ There
exists a }G\text{-equivariant map }F:V_{2}(\mathbb{R}^{3})\rightarrow
M(\alpha )\}.
\end{equation*}
\end{problem}

\noindent The well known \textquotedblleft extension of
scalars\textquotedblright\ equivalence from homological algebra, \cite{Brown}
Section III.3, has a very useful analogue in the category of $G$-spaces and $%
G$-equivariant maps, see the Section \ref{sec:geometry}. This equivalence
permits us in some cases to change the group and replace the original $G$%
-space by a new, more tractable topological space, see Proposition \ref%
{PropositionScalar}. Here we replace the dihedral group $D_{2n}$ by the
generalized quaternion group $Q_{4n}$ and the original configuration space $%
V_{2}(\mathbb{R}^{3})$ by the sphere $S^{3}.$

\noindent Let $S^{3}=S(\mathbb{H})=Sp(1)$ be the group of all unit
quaternions and let $\epsilon =\epsilon _{2n}=\cos \frac{\pi }{n}+i\sin
\frac{\pi }{n}\in S(\mathbb{H})$ be a root of unity. Group $\langle \epsilon
\rangle $ is a subgroup of $S(\mathbb{H})$ of the order $2n$. Then, the
\textit{generalized quaternion group, }\cite{CaEi}\textit{\ }p. 253\textit{,
}is the subgroup
\begin{equation*}
Q_{4n}=\{1,\epsilon ,\ldots ,\epsilon ^{2n-1},j,\epsilon j,\ldots ,\epsilon
^{2n-1}j\}
\end{equation*}%
of $S^{3}$ of order $4n$. Let $H=\{1,\epsilon ^{n}\}=\{1,-1\}\subset Q_{4n}$%
. Then, it is not hard to prove that the quotient group $Q_{4n}/H$ is
isomorphic to the dihedral group $D_{2n}$ of the order $2n$.

\begin{proposition}
\label{cor:korolar} There exists a $D_{2n}$-map $G:V_{2}(\mathbb{R}%
^{3})\rightarrow M(\alpha )$ if and only if there exists a $Q_{4n}$-map $%
F:S^{3}\rightarrow M(\alpha )$ where the generalized quaternion group\textit{%
\ }$Q_{4n}$ acts on $S^{3}$ as a subgroup of $Sp(1)\cong S^{3}$, and its
action on $W_{n}$ is described by $\epsilon (x_{1},\ldots
,x_{n})=(x_{2},\ldots ,x_{n},x_{1})$ and $j(x_{1},\ldots
,x_{n})=(x_{n},\ldots ,x_{2},x_{1})$.
\end{proposition}

\begin{proof}
Note that $\epsilon ^{n}=j^{2}$ acts trivially on $W_{n}$. Moreover, there
is an isomorphism $S^{3}/H\cong \mathbb{R}P^{3}\cong SO(3)\cong V_{2}(%
\mathbb{R}^{3})$. Finally, the $Q_{4n}/H\cong D_{2n}$ action on $%
S^{3}/H\cong V_{2}(\mathbb{R}^{3})$ coincides with the $D_{2n}$ action on $%
V_{2}(\mathbb{R}^{3})$ described in preceding section. Thus this proposition
is a direct consequence of the Proposition~\ref{PropositionScalar}, Section~%
\ref{sec:modifying}.
\end{proof}

\begin{remark}
Note that the $Q_{4n}$-action on $S^{3}$ is free. The $Q_{4n}$-action on $%
W_{n}$ is the restriction of the action on $\mathbb{R}^{n}$ described by
\begin{equation*}
\epsilon \cdot e_{i}=e_{(i+1)\func{mod}n}\text{ and }j\cdot e_{i}=e_{n-i+1}%
\text{.}
\end{equation*}
where $e_{1},..,e_{n}$ be the standard orthonormal basis in $\mathbb{R}^{n}$.
\end{remark}

In light of Proposition~\ref{cor:korolar}, Problem~\ref{prob:dihedral} is
equivalent to the following problem.

\begin{problem}
Describe the set
\begin{eqnarray*}
\Gamma _{G} =\{\alpha \in \tfrac{1}{n}\,\mathbb{N}^{4}\,|\,\text{ There
exists a }G\text{-equivariant map }F:S^{3}\rightarrow M(\alpha )\}.
\end{eqnarray*}
\end{problem}

\section{{Topological preliminaries}}

\noindent

\subsection{Obstruction Theory}

\label{sec:obs-theo}

\noindent Let $G$ be a subgroup of the generalized quaternion group $Q_{4n}$%
. A fundamental problem is to decide if there exists a $G$-equivariant map $%
F:S^{3}\rightarrow M(\alpha )$. The equivariant obstruction theory, as
presented in \cite{Dieck87} (see also \cite{guide2}, Sections 4.1--4.3), is
a versatile tool for studying this question. Recall that $[X,Y]$ is the set
of all homotopy classes of maps $f:X\longrightarrow Y$ while in case both $X$
and $Y$ are $G$-spaces, $[X,Y]_{G}$ is the corresponding set of all $G$%
-homotopy classes of $G$-equivariant maps.

\noindent Since the space $M=M(\alpha )$ is $1$-connected and consequently $%
2 $-simple in the sense that $\pi _{1}(M)$ acts trivially on $\pi _{2}(M)$,
then by \cite{Dieck87}, SectionII.3, and by Hurewicz theorem, $\pi
_{2}(M)\cong \lbrack S^{2},M]\cong H_{2}(M;\mathbb{Z})$ as $Q_{4n}$-modules.
Thus, the relevant part of the obstruction exact sequence, \cite{Dieck87},
\cite{guide2}, has the following form,

\begin{equation*}
\begin{array}{ccccc}
\lbrack S^{3},M]_{G} & \overset{\theta }{\longrightarrow } & \mathrm{Im}%
\left\{ [S_{(2)}^{3},M]_{G}\longrightarrow \lbrack S_{(1)}^{3},M]_{G}\right\}
& \overset{\tau }{\longrightarrow } & H_{G}^{3}(S^{3},H_{2}(M;\mathbb{Z})).%
\end{array}%
\end{equation*}%
Here $S_{(1)}^{3}$ and $S_{(2)}^{3}$ are respectively the $1$ and $2$%
-skeleton of the sphere $S^{3}$, relative to some
$Q_{4n}$-invariant simplicial or CW-structure. Our initial choice
is the simplicial structure arising from the join decomposition
$S^{3}=P_{2n}^{(1)}\ast P_{2n}^{(2)}$ where both $P_{2n}^{(1)}$
and $P_{2n}^{(2)}$ are regular $(2n)$-sided polygons, see
Figure~\ref{fig:taraba2}. More precisely, the vertices of
these polygons are respectively, $v_{i}:=\epsilon ^{i}a$ and $%
w_{i}:=\epsilon ^{i}ja,\,i=0,1,\ldots ,2n-1$, where $a$ is a fixed complex
number, $|a|=1$. Since $[S_{(1)}^{3},M]_{G}=\{\ast \}$ is a one-element set
and $[S_{(2)}^{3},M]_{G}$ is nonempty, the obstruction exact sequence
reduces to,
\begin{equation*}
\begin{array}{ccccc}
\lbrack S^{3},M]_{G} & \longrightarrow & \{\ast \} & \overset{\tau }{%
\longrightarrow } & H_{G}^{3}(S^{3},H_{2}(M;\mathbb{Z})).%
\end{array}%
\end{equation*}%
The exactness of this sequence means that the set $[S^{3},M]_{G}\neq
\emptyset $ if and only if a special element $\tau (\ast )\in
H_{G}^{3}(S^{3},H_{2}(M;\mathbb{Z}))$ is equal to zero. Note that $\tau
(\ast )$ depends only on $M$.

\noindent The class $\tau (\ast )$ can be evaluated by studying the
\textquotedblleft singular set\textquotedblright\ of a general position
equivariant map, see \cite{guide2}, Sections 4.1-4.3, for a brief overview.
A $G$-simplicial map $h:S^{3}\rightarrow W_{n}$ satisfies a
\textquotedblleft general position condition\textquotedblright\ if for each
simplex $\sigma $ in $S^{3}$ and any linear space $U$ in the arrangement $%
\mathcal{A}(\alpha )$
\begin{equation*}
h(\sigma )\cap U\neq \emptyset \Rightarrow \mathrm{dim}(\sigma )=\mathrm{dim}%
(h(\sigma ))=3,\mathrm{dim}(U)=n-3,\,h(\sigma )\cap U=\{pt\}\subset \mathrm{%
int(}h(\sigma )).
\end{equation*}%
It is not difficult to check that a generic map is in the general position
and that every equivariant, simplicial map can be put in general position by
a small perturbation. For each $G$-map $h:S^{3}\rightarrow W_{n\text{ \ }}$%
in the general position, there is an associated obstruction cocycle
\begin{equation*}
c(h)\in C_{G}^{3}(S^{3},{A})=\mathrm{Hom}_{G}(C_{3}(S^{3}),A),\,[c(h)]=\tau
(\ast ).
\end{equation*}%
If $\sigma $ is an oriented $3$-simplex in $S^{3}$, then $c(h)(\sigma )\in
H_{2}(M;\mathbb{Z})$ is the image $h_{\ast }([\partial (\sigma )])$ of the
fundamental class of $\partial (\sigma )\cong S^{2}$ by the map $h_{\ast
}:H_{2}(\partial (\sigma );\mathbb{Z})\rightarrow H_{2}(M;\mathbb{Z})$. To
find explicit form of the obstruction cocycle $c(h)\in
C_{G}^{3}(S^{3},H_{2}(M;\mathbb{Z}))=\mathrm{Hom}_{G}(C_{3}(S^{3}),H_{2}(M;%
\mathbb{Z}))$ we recall that by definition
\begin{equation*}
c(h)(\theta )\neq 0\,\Longleftrightarrow \,h(\theta )\cap (\cup \mathcal{A}%
(\alpha ))\neq \emptyset
\end{equation*}%
and $c(h)(\theta )=h_{\ast }([\partial \theta ])\in H_{2}(M;{\mathbb{Z}})$.

\noindent Before we go further, let us record for the future reference an
important property of the obstruction cocycle. More details on the
restriction and the transfer map can be found in \cite{Brown} Section~III.9.

\begin{proposition}
\label{Th.ObCoCycleTorzioni} Let $G$ be subgroup of $Q_{4n}$. The cohomology
class of the obstruction cocycle $c(h)$ is a torsion element of the group $%
H_{G}^{3}(S^{3},H_{2}(M;\mathbb{Z}))$.
\end{proposition}

\begin{proof}
Let $H$ be a subgroup of $G$. There exists a natural \textquotedblleft
restriction\textquotedblright\ map $r:H_{G}^{3}(S^{3},H_{2}(M;\mathbb{Z}%
))\rightarrow H_{H}^{3}(S^{3},H_{2}(M;\mathbb{Z})),$ which on the cochain
level is just the \textquotedblleft forgetful map" sending a $G$-cochain $%
c\in C_{G}^{3}(S^{3},H_{2}(M;\mathbb{Z}))$ to the same cochain interpreted
as a $H$-cochain. It follows from the definition of the obstruction cocycle
that $r(c(h))$ is the obstruction cocycle for the extension of a general
position $H$-map $h$. Moreover there exists a natural map $\tau
:H_{H}^{3}(S^{3},H_{2}(M;\mathbb{Z}))\rightarrow H_{G}^{3}(S^{3},H_{2}(M;%
\mathbb{Z}))$ in the opposite direction called the transfer map. It is
known, \cite{Brown} Section III.9. Proposition 9.5.(ii)), that the
composition of the restriction with the transfer is just multiplication by
the index $[G:H]$:
\begin{equation*}
\begin{array}{ccccc}
H_{G}^{3}(S^{3},H_{2}(M;\mathbb{Z})) & \rightarrow & H_{H}^{3}(S^{3},H_{2}(M;%
\mathbb{Z})) & \rightarrow & H_{G}^{3}(S^{3},H_{2}(M;\mathbb{Z})) \\
\lbrack c(h)_{G}] & \longmapsto & [c(h)_{H}] & \longmapsto & [G:H]\cdot
\lbrack c(h)_{G}]\text{.}%
\end{array}%
\end{equation*}%
Note that if $H$ is the trivial group, the cohomology class of the
obstruction cocycle $[c(h)_{H}]$ is zero. This implies that $[G:H]\cdot
\lbrack c(h)_{G}]=0$ in $H_{G}^{3}(S^{3},H_{2}(M;\mathbb{Z}))$, i.e. $%
[c(h)_{G}]$ is a torsion element of the group $H_{G}^{3}(S^{3},H_{2}(M;%
\mathbb{Z}))$.
\end{proof}

\begin{remark}
The previous result explains why we pay a special attention in subsequent
sections to the torsion part of the group $H_{G}^{3}(S^{3},H_{2}(M;\mathbb{Z}%
))$.
\end{remark}

\subsection{The Z-\v{Z} formula and point classes of $M(\mathcal{A}(\protect%
\alpha ))$}

\noindent Before we continue with the calculation of the obstruction cocycle
$c_{G}(h)\in C_{G}^{3}(S^{3},H_{2}(M;\mathbb{Z}))$, we should gather
together more information about the $G$-module $H_{2}(M;\mathbb{Z})$. Since
every subspace $K\in \mathcal{A}(\alpha )$ has the codimension in $W_{n}$
greater or equal to $3$, the fundamental group $\pi _{1}(M)$ of the
complement $M=W_{n}\setminus \cup \mathcal{A}(\alpha )$ is trivial. Thus, $%
H_{1}(M;\mathbb{Z})=0$ and the first nontrivial homology group of the
complement is $H_{2}(M;\mathbb{Z})$.

\bigskip

\noindent Before we begin the analysis of the homology group $H_{2}(M;%
\mathbb{Z})$of the complement $M=W_{n}\setminus \cup \mathcal{A}(\alpha )$,
we make a few general observations about representations of homology classes
of complements of arrangements. Suppose that $L_{1},L_{2},\ \ldots ,L_{k}$
is a collection of $a$-dimensional linear subspaces in an $(a+b)$%
-dimensional, Euclidean space $V$. Let $\mathcal{A}$ be the associated
arrangement with the intersection poset $P_{\mathcal{A}}$, compactified
union (link) $\hat{D}(\mathcal{A})=\cup \mathcal{A}\cup \{+\infty \}\subset
V\cup \{+\infty \}\cong S^{a+b}$, and the complement $M(\mathcal{A}%
)=V\setminus \cup \mathcal{A}$. For a given point $x\in L_{i}\setminus \cup
_{j\neq i}L_{j}$, let $D_{\epsilon }(x)=x+D_{\epsilon }$ be a small disc
around $x$, where $D_{\epsilon }=\{y\in L_{i}^{\perp }\mid \langle
y,y\rangle \leq \epsilon \}$. \textquotedblleft Small\textquotedblright\
means that $\epsilon $ is chosen so that $D_{\epsilon }(x)\cap (\cup _{j\neq
i}L_{j})=\emptyset $. We assume that $D_{\epsilon }(x)$ is oriented,
typically by the orientation inherited from some orientations on $V$ and $%
L_{i}$ which are prescribed in advance. The fundamental class of the pair $%
(D_{\epsilon }(x),\partial D_{\epsilon }(x))$ determines a homology class in
$H_{b}(V,M(\mathcal{A});\mathbb{Z})$, which we denote by $[x]$ and call the%
\textit{\ point class} of $x$. Note that by the Excision axiom, $[x]$ does
not depend on $\epsilon $. Moreover, by the Homotopy axiom, $[x]$ does not
change if $x$ is moved inside a connected component of $L_{i}\setminus \cup
_{j\neq i}L_{j}$. Similarly, in light of the isomorphism $H_{b}(V,M(\mathcal{%
A}))\rightarrow H_{b-1}(M(\mathcal{A}))$, the class ${[[}x{]]}:=\partial
\lbrack x]$, which is also called the \textit{point class }of\emph{\ }$x$,
has all these invariance properties as well.

\noindent Let us show that the class $[x]$ is always nontrivial. By the
Ziegler-\v{Z}ivaljevi\'{c} formula \cite{ZZ}, the homotopy type of the
one-point compactification $\hat{D}(\mathcal{A})=\cup \mathcal{A}\cup
\{+\infty \}$ has the wedge decomposition of the form,
\begin{equation*}
\hat{D}(\mathcal{A})\simeq \hat{L}_{1}\vee \hat{L}_{2}\vee \ldots \vee \hat{L%
}_{k}\vee \ldots
\end{equation*}%
where the displayed factors correspond to elements $p\in P(\mathcal{A})$ of
the minimum dimension $d(p)=a$. Let $c(\hat{L}_{i})\in H^{b}(V,M(\mathcal{A}%
);\mathbb{Z})$ be the cohomology class which is Poincar\'{e}-Alexander dual
of the fundamental homology class $[\hat{L}_{i}]\in H_{a}(\hat{D}(\mathcal{A}%
);\mathbb{Z})$ associated with the sphere $\hat{L}_{i}$ in $\hat{V}$. Then $%
c(\hat{L}_{i})([N])\in \mathbb{Z}$ is essentially the intersection number $[%
\hat{l}_{i}]\cap \lbrack N]$, whenever this number is correctly defined, for
example if $N$ is a manifold and the intersection $\hat{L}_{i}\cap N$ is
transversal. From here we see that $c(\hat{L}_{i})([x])=1$, hence $[x]$ must
be nontrivial. All this also applies to the class ${[[}x]]=\partial \lbrack
x]\in H_{b-1}(M(\mathcal{A}))$. Consequently we have the following
proposition.

\begin{proposition}
\label{lem:main2} If the codimension of $L_{ij}:=L_{i}\cap L_{j}$ in $L_{i}$
is greater than $1$ for each $j\neq i$ and $x\in L_{i}\setminus \cup _{j\neq
i}L_{j}$, then the point class $[x]$ is a well defined class in $H_{b}(V,M(%
\mathcal{A}))$ which does not depend on $x\in L_{i}$ whatsoever.
\end{proposition}

\noindent What happens when $\mathrm{codim}_{L_{i}}(L_{i}\cap L_{j})=1$ for
some $j\neq i$? Then $L_{ij}$, being a hyperplane in both $L_{i}$ and $L_{j}$%
, decomposes these spaces into the union of closed halfspaces, $%
L_{i}=L_{i}^{1}\cup L_{i}^{2}$ and $L_{j}=L_{j}^{1}\cup L_{j}^{2}$
respectively. By gluing the halfspaces $L_{i}^{1}$ and $L_{j}^{1}$ along the
common boundary $L_{ij}$ and by adding the infinite point $+\infty $, one
obtains a sphere $K\subset \hat{V}$. Recall that the decomposition in
Ziegler-\v{Z}ivaljevi\'{c} formula \cite{ZZ} involves a choice of generic
points in all elements of the arrangement $\mathcal{A}$. The points can be
chosen in halfspaces $L_{i}^{1}$ and $L_{j}^{1}$, which implies that the
sphere $K$ appears in the decomposition as the factor $\hat{L}_{ij}\ast
\Delta (P_{<p_{ij}})\cong S^{a-1}\ast S^{0}\cong S^{a}$, where $p_{ij}\in P$
corresponds to $L_{ij}\in \mathcal{A}$. This guarantees that the fundamental
class $[K]$ of $K$ is nontrivial in $H_{a}(\hat{D}(\mathcal{A}))$. Let $%
c(K)\in H^{b}(V,M(\mathcal{A});\mathbb{Z})$ be the Poincar\'{e}-Alexander
dual to $[K]$. Then $c(K)([x_{1}])=\pm 1$ if $x_{1}$ is in the interior of
the halfspace $L_{i}^{1}$ but $c(K)([x_{2}])=0$ if $x_{2}$ belongs to the
interior of the complementary halfspace $L_{i}^{2}$. This observation leads
to the following proposition.

\begin{proposition}
Let $\mathrm{codim}_{L_{i}}(L_{i}\cap L_{j})=1$ for some $j\neq i$ where $%
L_{i}=L_{i}^{1}\cup L_{i}^{2}$, $L_{j}=L_{j}^{1}\cup L_{j}^{2}$ are
decompositions in closed halfspaces. If points $x_{1}$ and $x_{2}$ belong to
interiors of complementary halfspaces $L_{i}^{1}$ and $L_{i}^{2}$then $%
[x_{1}]\neq \lbrack x_{2}]$.
\end{proposition}

\section{{Calculation of the obstruction cocycles}}

\noindent Now we are ready to compute the homology class of the cocycle $%
c_{G}(h)\in C_{G}^{3}(S^{3},H_{2}(M;\mathbb{Z}))$. As before $G$ is a
subgroup of the generalized quaternion group. Our primary interest is in the
complete group $Q_{4n}$ and the cyclic subgroup $\mathbb{Z}_{n}=\{1,\epsilon
^{2},\ldots ,\epsilon ^{2n-2}\}$. In the latter case we complete the
calculations started in \cite{VreZiv2002}.

\subsection{The obstruction cocycle}

\subsubsection{General position $Q_{4n}$-maps and their singular sets}

\label{sec:genpos}

\noindent Let us start with the description of a general position,
simplicial $Q_{4n}$-map $h:S^{3}\rightarrow W_{n}$, where the sphere has the
simplicial structure $S^{3}=P_{2n}^{(1)}\ast P_{2n}^{(2)}$, described in
Section~\ref{sec:obs-theo} and depicted in Figure~\ref{fig:taraba2}. Let $%
e_{1},..,e_{n}$ be the standard orthonormal basis in $\mathbb{R}^{n}$ and
let $x_{1},..,x_{n}$ be associated dual linear functions. Let $%
\{u_{1},..,u_{n}\}$, where $u_{i}=e_{i}-e$ and $e=\frac{1}{n}%
\sum_{r=1}^{n}e_{r}$, be the vertex set of a regular simplex $\Delta _{n-1}$
in $W_{n}$. If the map $h$ is prescribed in advance on the vertex $a\in
P_{2n}^{(1)}\ast P_{2n}^{(2)}$, say if $h(a)=u_{1}$, and if we require that
it is a simplicial $Q_{4n}$-map, then everything else is completely
determined. For example,
\begin{eqnarray*}
h(\epsilon ^{i}a) &=&\epsilon ^{i}\cdot h(a)=\epsilon ^{i}\cdot
u_{1}=\epsilon ^{i}\cdot (e_{1}-e)=e_{i\func{mod}n+1}-e=u_{i\func{mod}n+1} \\
h(ja) &=&j\cdot h(a)=j\cdot u_{1}=j\cdot (e_{1}-e)=e_{n}-e=u_{n} \\
h(\epsilon ^{i}ja) &=&\epsilon ^{i}j\cdot h(a)=\epsilon ^{i}\cdot
(e_{n}-e)=e_{(i+n)\func{mod}n}-e=e_{i}-e=u_{i\func{mod}n}.
\end{eqnarray*}%
To see that $h$ is a general position map, it is sufficient to test for
which simplexes $\sigma =\sigma _{1}\ast \sigma _{2}\subset P_{2n}\ast
P_{2n} $ the image $h(\sigma )$ intersects the subspace $L=L(\alpha )$.
First we observe that $h(P_{2n}\ast P_{2n})\subseteq \mathrm{sk}_{3}(\Delta
_{n-1})$ and that $L\cap h(P_{2n}\ast P_{2n})=\{y_{1},y_{2}\}$, where
\begin{eqnarray*}
y_{1} &=&\frac{a_{1}}{n}u_{a_{1}}+\frac{a_{2}}{n}u_{a_{1}+1}+\frac{a_{3}}{n}%
u_{a_{1}+a_{2}+a_{3}}+\frac{a_{4}}{n}u_{a_{1}+a_{2}+a_{3}+1} \\
y_{2} &=&\frac{a_{2}}{n}u_{a_{1}+a_{2}}+\frac{a_{3}}{n}u_{a_{1}+a_{2}+1}+%
\frac{a_{4}}{n}u_{n}+\frac{a_{1}}{n}u_{1}.
\end{eqnarray*}%
Thus, there are only two $3$-simplices
\begin{equation*}
\tau
_{1}=[u_{a_{1}},u_{a_{1}+1};u_{a_{1}+a_{2}+a_{3}},u_{a_{1}+a_{2}+a_{3}+1}]%
\text{ and }\tau _{2}=[u_{a_{1}+a_{2}},u_{a_{1}+a_{2}+1};u_{n},u_{1}]
\end{equation*}%
in the simplex $\Delta _{n-1}$ which intersect $L$. To find the singular set
$h^{-1}(\cup \mathcal{A}(\alpha ))\subseteq $ $S^{3}$, we have to detect all
$3$-simplices $\sigma =\sigma _{1}\ast \sigma _{2}\subset P_{2n}\ast P_{2n}$
in the sphere $S^{3}$ with the property $\sigma \cap h^{-1}(L)\neq \emptyset
$, or in other words simplices $\sigma $ such that either $h(\sigma )=\tau
_{1}$ or $h(\sigma )=\tau _{2}$. This leads to the following systems of
equations:
\begin{equation*}
\begin{array}{lll}
h(\epsilon ^{i}a)=u_{i\func{mod}n+1}=u_{a_{1}} & \text{ and } & h(\epsilon
^{i}ja)=u_{i\func{mod}n}=u_{a_{1}+a_{2}+a_{3}} \\
h(\epsilon ^{i}a)=u_{i\func{mod}n+1}=u_{a_{1}+a_{2}+a_{3}} & \text{ and } &
h(\epsilon ^{i}ja)=u_{i\func{mod}n}=u_{a_{1}} \\
h(\epsilon ^{i}a)=u_{i\func{mod}n+1}=u_{a_{1}+a_{2}} & \text{ and } &
h(\epsilon ^{i}ja)=u_{i\func{mod}n}=u_{n} \\
h(\epsilon ^{i}a)=u_{i\func{mod}n+1}=u_{n} & \text{ and } & h(\epsilon
^{i}ja)=u_{i\func{mod}n}=u_{a_{1}+a_{2}}.%
\end{array}%
\end{equation*}%
There are 16, not necessarily different, $3$-simplices $\theta _{1,}\ldots
,\theta _{16}$ in $P_{2n}\ast P_{2n}$ which nontrivially intersect the
singular subset $h^{-1}(L)$. Here is a complete list where $v_{i}=\epsilon
^{i}a$, $w_{i}=\epsilon ^{i}ja$ for $i\in \{0,..,2n-1\}$ and $P=a_{1}$, $%
Q=a_{1}+a_{2}$, $R=a_{1}+a_{2}+a_{3}$.

\begin{center}
$%
\begin{array}{ll}
\theta _{1}=[v_{P-1},v_{P};w_{R},w_{R+1}], & \theta
_{2}=[v_{n+P-1},v_{n+P};w_{R},w_{R+1}], \\
\theta _{3}=[v_{P-1},v_{a_{1}};w_{n+R},w_{n+R+1}], & \theta
_{4}=[v_{n+P-1},v_{n+P};w_{n+R},w_{n+R+1}], \\
\theta _{5}=[v_{R-1},v_{R};w_{P},w_{P+1}], & \theta
_{6}=[v_{n+R-1},v_{n+R};w_{P},w_{P+1}], \\
\theta _{7}=[v_{R-1},v_{R};w_{n+P},w_{n+P+1}], & \theta
_{8}=[v_{n+R-1},v_{n+R};w_{n+P},w_{n+P+1}], \\
\theta _{9}=[v_{Q-1},v_{Q};w_{n},w_{n+1}], & \theta
_{10}=[v_{n+Q-1},v_{n+Q};w_{n},w_{n+1}], \\
\theta _{11}=[v_{Q-1},v_{Q};w_{0},w_{1}], & \theta
_{12}=[v_{n+Q-1},v_{n+Q};w_{0},w_{1}], \\
\theta _{13}=[v_{n-1},v_{n};w_{Q},w_{Q+1}], & \theta
_{14}=[v_{2n-1},v_{0};w_{Q},w_{Q+1}], \\
\theta _{15}=[v_{n-1},v_{n};w_{n+Q},w_{n+Q+1}], & \theta
_{16}=[v_{2n-1},v_{0};w_{n+Q},w_{n+Q+1}].%
\end{array}%
$
\end{center}

Note that
\begin{equation*}
h(\theta _{i})=\tau _{1}\text{ if }i\in \{1,..,8\}\text{ \ and }h(\theta
_{i})=\tau _{2}\text{ if }i\in \{9,..,16\}\text{.}
\end{equation*}

\subsubsection{The obstruction cocycle $c_{Q_{4n}}(h)$}

\label{sec:obscoc}

\noindent We have made all necessary preparations and now we are finally
ready to compute the obstruction cocycle $c_{Q_{4n}}(h)$. For a $3$-simplex $%
\theta $ in $S^{3}$,
\begin{equation*}
c_{Q_{4n}}(h)(\theta )=\sum_{y\in h(\theta )\cap \cup \mathcal{A}(\alpha )}%
\mathrm{I}(h(\theta ),L_{y})[y]
\end{equation*}%
where $\mathrm{I}(h(\theta ),L_{y})$ is the intersection number of the
oriented simplex $h(\theta )$ and the appropriate oriented element $L_{y}$ ($%
L_{y}\cap h(\theta )=\{y\}$) of the arrangement $\mathcal{A}(\alpha )$ and $%
[y]\in H_{2}(M;\mathbb{Z})$ is a point class.

\noindent In order to simplify the computation it is convenient to replace
the simplicial complex structure $P_{2n}\ast P_{2n}$ on $S^{3}$ and its
chain complex $\{C_{i}(S^{3},\mathbb{Z})\}_{i=0}^{3}$\ by a much more
economical $Q_{4n}$-invariant, cell complex structure on $S^{3}$ with the
chain complex $\{D_{i}(S^{3},\mathbb{Z})\}_{i=0}^{3}$, described in the
Appendix, Section{\ \ref{sec:geometry}}. The obstruction class computed
relative to the new cell complex structure is denoted by $c_{Q_{4n}}^{\prime
}(h)$.

\noindent To evaluate the obstruction cocycle $c_{Q_{4n}}^{\prime }(h)$ on
the cell $e$, keeping in mind natural maps of chain complexes $\{C_{i}(S^{3},%
\mathbb{Z})\}_{i=0}^{3}$ and $\{D_{i}(S^{3},\mathbb{Z})\}_{i=0}^{3}$, we
have to find how many simplexes of the form $g\cdot \theta _{i}$, where $%
g\in Q_{4n}$ and $i\in \{1,..,16\}$, belong to the cell $e$. It is found by
inspection that

\begin{eqnarray*}
\sigma _{1} &=&[v_{a_{1}+a_{4}-1},v_{a_{1}+a_{4}};w_{0},w_{1}]=\epsilon
^{-a_{1}-a_{2}-a_{3}}\theta _{2}=j\epsilon ^{-a_{1}}\theta _{1}, \\
\sigma _{2} &=&[v_{a_{2}+a_{3}-1},v_{a_{2}+a_{3}};w_{0},w_{1}]=\epsilon
^{-a_{1}}\theta _{5}=j\epsilon ^{-a_{1}-a_{2}-a_{3}}\theta _{7}, \\
\sigma _{3} &=&[v_{a_{2}+a_{1}-1},v_{a_{2}+a_{1}};w_{0},w_{1}]=\theta
_{11}=j\epsilon ^{-a_{1}-a_{2}}\theta _{9}, \\
\sigma _{4} &=&[v_{a_{3}+a_{4}-1},v_{a_{3}+a_{4}};w_{0},w_{1}]=\epsilon
^{-a_{1}-a_{2}}\theta _{13}=j\theta _{14},
\end{eqnarray*}%
is the complete list of these simplices. In other words there are four (not
necessarily different \textbf{!}) simplices $\sigma _{1}$, $\sigma _{2}$, $%
\sigma _{3}$, $\sigma _{4}$ in the cell $e$ such that $h(\sigma _{i})\cap
D(\alpha )\neq \varnothing $ for each $i\in \{1,2,3,4\}$. Also observe that
for each of the simplices $\sigma _{1}$, $\sigma _{2}$, $\sigma _{3}$, $%
\sigma _{4}$ the intersection $h(\sigma _{i})\cap $ $D(\alpha )$ has two
(not necessarily different \textbf{!}) points. More explicitly these points
are%
\begin{eqnarray*}
x_{11} &=&\tfrac{a_{1}}{n}v_{a_{1}+a_{4}-1}+\tfrac{a_{2}}{n}v_{a_{1}+a_{4}}+%
\tfrac{a_{3}}{n}w_{0}+\tfrac{a_{4}}{n}w_{1}\in \sigma _{1}\cap
h^{-1}(\epsilon ^{-a_{1}-a_{2}-a_{3}}L), \\
x_{12} &=&\tfrac{a_{4}}{n}v_{a_{1}+a_{4}-1}+\tfrac{a_{3}}{n}v_{a_{1}+a_{4}}a+%
\tfrac{a_{2}}{n}w_{0}+\tfrac{a_{1}}{n}w_{1}\in \sigma _{1}\cap
h^{-1}(j\epsilon ^{-a_{1}}L), \\
x_{21} &=&\tfrac{a_{3}}{n}v_{a_{2}+a_{3}-1}+\tfrac{a_{4}}{n}v_{a_{2}+a_{3}}+%
\tfrac{a_{1}}{n}w_{0}+\tfrac{a_{2}}{n}w_{1}\in \sigma _{2}\cap
h^{-1}(\epsilon ^{-a_{1}}L), \\
x_{22} &=&\tfrac{a_{2}}{n}v_{a_{2}+a_{3}-1}+\tfrac{a_{1}}{n}v_{a_{2}+a_{3}}+%
\tfrac{a_{4}}{n}w_{0}+\tfrac{a_{3}}{n}w_{1}\in \sigma _{2}\cap
h^{-1}(j\epsilon ^{-a_{1}-a_{2}-a_{3}}L), \\
x_{31} &=&\tfrac{a_{2}}{n}v_{a_{2}+a_{1}-1}+\tfrac{a_{3}}{n}v_{a_{2}+a_{1}}+%
\tfrac{a_{4}}{n}w_{0}+\tfrac{a_{1}}{n}w_{1}\in \sigma _{3}\cap h^{-1}(L), \\
x_{32} &=&\tfrac{a_{1}}{n}v_{a_{2}+a_{1}-1}+\tfrac{a_{4}}{n}v_{a_{2}+a_{1}}+%
\tfrac{a_{3}}{n}w_{0}+\tfrac{a_{2}}{n}w_{1}\in \sigma _{3}\cap
h^{-1}(j\epsilon ^{-a_{1}-a_{2}}L), \\
x_{41} &=&\tfrac{a_{4}}{n}v_{a_{3}+a_{4}-1}+\tfrac{a_{1}}{n}v_{a_{3}+a_{4}}+%
\tfrac{a_{2}}{n}w_{0}+\tfrac{a_{3}}{n}w_{1}\in \sigma _{4}\cap
h^{-1}(\epsilon ^{-a_{1}-a_{2}}L), \\
x_{42} &=&\tfrac{a_{3}}{n}v_{a_{3}+a_{4}-1}+\tfrac{a_{2}}{n}v_{a_{3}+a_{4}}+%
\tfrac{a_{1}}{n}w_{0}+\tfrac{a_{4}}{n}w_{1}\in \sigma _{4}\cap h^{-1}(jL),
\end{eqnarray*}

\noindent then
\begin{eqnarray*}
h(\sigma _{1})\cap D(\alpha ) &=&(h(\sigma _{1})\cap \epsilon
^{-a_{1}-a_{2}-a_{3}}L)\cup (h(\sigma _{1})\cap j\epsilon
^{-a_{1}}L)=\{h(x_{11}),h(x_{12})\}; \\
h(\sigma _{2})\cap D(\alpha ) &=&(h(\sigma _{2})\cap \epsilon
^{-a_{1}}L)\cup (h(\sigma _{2})\cap j\epsilon
^{-a_{1}-a_{2}-a_{3}}L)=\{h(x_{21}),h(x_{22})\}; \\
h(\sigma _{3})\cap D(\alpha ) &=&(h(\sigma _{3})\cap L)\cup (h(\sigma
_{3})\cap j\epsilon ^{-a_{1}-a_{2}}L)=\{h(x_{31}),h(x_{32})\}; \\
h(\sigma _{4})\cap D(\alpha ) &=&(h(\sigma _{4})\cap \epsilon
^{-a_{1}-a_{2}}L)\cup (h(\sigma _{4})\cap jL)=\{h(x_{41}),h(x_{42})\}
\end{eqnarray*}%
and
\begin{equation*}
h(e)\cap D(\alpha
)=%
\{h(x_{11}),h(x_{12}),h(x_{21}),h(x_{22}),h(x_{31}),h(x_{32}),h(x_{41}),h(x_{42})\}.
\end{equation*}%
Now we focus on the cardinality of the intersection $h(e)\cap D(\alpha ).$As
it turns out, the spaces $D(\alpha )$ and the corresponding cardinalities
are very sensitive to the change of integers $a_{1},a_{2},a_{3},a_{4.}$ We
pay a special attention to the question whether the test map $h$ is in
general position since this is a key technical assumption used in
computations. We discuss several separate cases.

\textbf{(A)} Let $a_{2}=a_{4}$ and $a_{1}\neq a_{3}$, in which case $\sigma
_{1}=\sigma _{3}$, $\sigma _{2}=\sigma _{4}$ and $x_{11}=x_{32}$, $%
x_{12}=x_{31}$, $x_{21}=x_{42}$, $x_{22}=x_{41}$. These identities imply
that
\begin{equation*}
\begin{array}{ll}
h(x_{11})=h(x_{32})\in \epsilon ^{-a_{1}-a_{2}-a_{3}}L\cap j\epsilon
^{-a_{1}-a_{2}}L\text{, } & h(x_{12})=h(x_{31})\in j\epsilon ^{-a_{1}}L\cap L%
\text{, } \\
h(x_{21})=h(x_{42})\in \epsilon ^{-a_{1}}L\cap jL & h(x_{22})=h(x_{41})\in
j\epsilon ^{-a_{1}-a_{2}-a_{3}}L\cap \epsilon ^{-a_{1}-a_{2}}L%
\end{array}%
\end{equation*}%
A question arises whether $h$ is in general position. However, cf. the proof
of Theorem \ref{Th.KoinvarijanteQ4n} (J), there is a relation $\epsilon
^{-a_{1}}L=jL$ which implies%
\begin{eqnarray*}
j\epsilon ^{-a_{1}-a_{2}}L &=&\epsilon ^{(n-1)(-a_{1}-a_{2})}jL=\epsilon
^{a_{1}+a_{2}}jL=\epsilon ^{a_{2}}L=\epsilon ^{a_{4}}L=\epsilon
^{-a_{1}-a_{2}-a_{3}}L\text{,} \\
\epsilon ^{-a_{1}}L &=&jL\text{ } \\
j\epsilon ^{-a_{1}-a_{2}-a_{3}}L &=&j\epsilon ^{a_{4}}L=\epsilon
^{(n-1)a_{4}}jL=\epsilon ^{-a_{4}-a_{1}}L=\epsilon ^{-a_{2}-a_{1}}L \\
j\epsilon ^{-a_{1}}L &=&\epsilon ^{(n-1)(-a_{1})}jL=\epsilon ^{a_{1}}jL=L
\end{eqnarray*}%
and shows that $h$ is indeed in general position.

\textbf{(B) }Let $a_{2}=a_{4}\neq a_{1}=a_{3}$, in which case $\sigma
_{1}=\sigma _{3}=\sigma _{2}=\sigma _{4}$ and $x_{11}=x_{32}=x_{21}=x_{42}$,
$x_{12}=x_{31}=x_{22}=x_{41}$. In this case%
\begin{eqnarray*}
h(x_{11}) &=&h(x_{32})=h(x_{21})=h(x_{42})\in \epsilon
^{-a_{1}-a_{2}-a_{3}}L\cap j\epsilon ^{-a_{1}-a_{2}}L\cap \epsilon
^{-a_{1}}L\cap jL, \\
h(x_{12}) &=&h(x_{31})=h(x_{22})=h(x_{41})\in j\epsilon ^{-a_{1}}L\cap L\cap
j\epsilon ^{-a_{1}-a_{2}-a_{3}}L\cap \epsilon ^{-a_{1}-a_{2}}L,
\end{eqnarray*}%
and again, there is a question of the genericity of $h$. The proof of
Theorem \ref{Th.KoinvarijanteQ4n}, (B) shows that there exist relations $%
\epsilon ^{a_{1}+a_{2}}L=L$ and $\epsilon ^{a_{2}}L=jL$ which again implies
that $\ h$ is in general position%
\begin{equation*}
\epsilon ^{-a_{1}}L=\epsilon ^{a_{2}}L=jL=j\epsilon
^{a_{1}+a_{2}}L=j\epsilon ^{-a_{1}-a_{2}}L=\epsilon ^{-a_{1}-a_{2}-a_{3}}L.
\end{equation*}

\textbf{(C) }Let $a_{1}=a_{3}$ and $a_{2}\neq a_{4}$, then $\sigma
_{1}=\sigma _{4}$, $\sigma _{2}=\sigma _{3}$ and $x_{11}=x_{42}$, $%
x_{12}=x_{41}$, $x_{21}=x_{32}$, $x_{22}=x_{31}$. Like in the case (A), it
is not difficult to prove that $h$ is again in general position.

\textbf{(D) }Let $a_{1}=a_{2}=a_{3}=a_{4}$, then $\sigma _{1}=\sigma
_{3}=\sigma _{2}=\sigma _{4}$ and $%
x_{11}=x_{32}=x_{21}=x_{42}=x_{12}=x_{31}=x_{22}=x_{41}$. Thus,%
\begin{equation*}
h(x_{11})=..=h(x_{41})\in \epsilon ^{-a_{1}-a_{2}-a_{3}}L\cap j\epsilon
^{-a_{1}-a_{2}}L\cap \epsilon ^{-a_{1}}L\cap jL\cap j\epsilon ^{-a_{1}}L\cap
L\cap j\epsilon ^{-a_{1}-a_{2}-a_{3}}L\cap \epsilon ^{-a_{1}-a_{2}}L
\end{equation*}%
and again Theorem \ref{Th.KoinvarijanteQ4n}, (A) helps with relations $%
\epsilon ^{a_{1}}L=L$ and $jL=L$. Therefore%
\begin{equation*}
j\epsilon ^{-a_{1}-a_{2}-a_{3}}L=j\epsilon ^{-a_{1}-a_{2}}L=j\epsilon
^{-a_{1}}L=jL=L=\epsilon ^{-a_{1}}L=\epsilon ^{-a_{1}-a_{2}-a_{3}}L=\epsilon
^{-a_{1}-a_{2}}L
\end{equation*}%
and $h$ is in general position.

\textbf{(E) }Let $a_{1}=a_{2}\neq a_{3}=a_{4}$, then $\sigma _{1}=\sigma
_{2} $ and $x_{11}=x_{22}$, $x_{12}=x_{21}$. In this case%
\begin{eqnarray*}
h(x_{11}) &=&h(x_{22})\in \epsilon ^{-a_{1}-a_{2}-a_{3}}L\cap j\epsilon
^{-a_{1}-a_{2}-a_{3}}L \\
h(x_{12}) &=&h(x_{21})\in j\epsilon ^{-a_{1}}L\cap \epsilon ^{-a_{1}}L.
\end{eqnarray*}%
With the relations $jL=\epsilon ^{2a_{3}}L$ (Theorem \ref%
{Th.KoinvarijanteQ4n}, (I)) we have%
\begin{eqnarray*}
j\epsilon ^{-a_{1}-a_{2}-a_{3}}L &=&j\epsilon ^{a_{4}}L=\epsilon
^{(n-1)a_{4}}jL=\epsilon ^{-a_{4}+2a_{3}}L=\epsilon ^{a_{4}}L=\epsilon
^{-a_{1}-a_{2}-a_{3}}L, \\
j\epsilon ^{-a_{1}} &=&\epsilon ^{(n-1)(-a_{1})}jL=\epsilon
^{a_{1}+2a_{3}}L=\epsilon ^{-a_{1}}L,
\end{eqnarray*}%
and conclude that $h$ is in general position again.

\textbf{(F)} Let $a_{2}=a_{3}\neq a_{1}=a_{4}$, then $\sigma _{3}=\sigma
_{4} $ and $x_{31}=x_{42}$, $x_{32}=x_{41}$. Like in the preceding case (E),
it is not hard to prove that $h$ is in general position.

\textbf{(G) }In all the remaining cases
\begin{equation*}
h(e)\cap D(\alpha
)=%
\{h(x_{11}),h(x_{12}),h(x_{21}),h(x_{22}),h(x_{31}),h(x_{32}),h(x_{41}),h(x_{42})\}.
\end{equation*}

\begin{theorem}
\label{Th.ObCycleQ4n}Let $c_{Q_{4n}}^{\prime }(h)\in
C_{Q_{4n}}^{3}(S^{3},H_{2}(M;\mathbb{Z}))$ be the (second) obstruction
cocycle for the map $h$ defined above. Let $K$ be the additive subgroup of $%
H_{2}(M;\mathbb{Z})$ generated by the elements of the form $g\cdot x-x$, $%
g\in Q_{4n}$, $x\in H_{2}(M;\mathbb{Z})$. As usual, the group of
coinvariants \ \cite{Brown} \ is $H_{2}(M;\mathbb{Z})_{Q_{4n}}=H_{2}(M;%
\mathbb{Z})/K$.

(A) Let $(a_{1},a_{2},a_{3},a_{4})=(p,p,p,p)$. Then $c_{Q_{4n}}^{\prime
}(h)(e)=[[y_{2}]]+K$ $\in H_{2}(M;\mathbb{Z})_{Q_{4n}}$.

(B) Let $(a_{1},a_{2},a_{3},a_{4})=(p,p,p,2p)$. Then $[c_{Q_{4n}}^{\prime
}(h)(e)]=2((1+(-1)^{{\binom{5p }{2}}}))[[y_{1}]]+K$ in $H_{2}(M;\mathbb{Z}%
)_{Q_{4n}}$.

(C) Let $(a_{1},a_{2},a_{3},a_{4})\neq (p,p,p,p)$. Then there exists $z\in
H_{2}(M;\mathbb{Z})$ such that $[c_{Q_{4n}}^{\prime }(h)(e)]=2z+K$ in $%
H_{2}(M;\mathbb{Z})_{Q_{4n}}$.
\end{theorem}

\begin{proof}
Before we calculate obstruction cocycle in all these cases recall that
\begin{equation*}
c_{Q_{4n}}(h)(\theta )=\sum_{y\in h(\theta )\cap D(\alpha )}\mathrm{I}%
(h(\theta ),L_{y})[y]
\end{equation*}%
where $\mathrm{I}(h(\theta ),L_{y})$ is the intersection number of the
oriented simplex $h(\theta )$ and the appropriate oriented element $L_{y}$ ($%
L_{y}\cap h(\theta )=\{y\}$) of the arrangement $\mathcal{A}(\alpha )$.
Also, for every $g\in Q_{4n}$ there is an identity
\begin{equation*}
\mathrm{I}(h(g\theta ),L_{g\cdot y})=\mathrm{\det }(g)\,\mathrm{I}(h(\theta
),L_{y})
\end{equation*}%
where $\det (\epsilon )=(-1)^{n+1}$, $\det (j)=(-1)^{\binom{(n}{2)}}$ and $%
\det (g\cdot g^{\prime })=\det (g)\cdot \det (g^{\prime })$. Following
calculations from \cite{VreZiv2002}, we have:
\begin{equation*}
\mathrm{I}(h(\theta _{i}),L)=\left\{
\begin{array}{ll}
{1} & ,{i\in \{1,..,8\}} \\
{-1} & ,{i\in \{9,..,16\}}%
\end{array}%
\right. ,\text{ where }h(\theta _{i})\cap L=h(x_{\theta _{i}})=\left\{
\begin{array}{ll}
{y}_{1} & ,\,{i\in \{1,..,8\}} \\
{y}_{2} & ,\,{i\in \{9,..,16\}}%
\end{array}%
\right. .
\end{equation*}%
From the first equation we get
\begin{equation*}
c_{Q_{4n}}^{\prime }(h)(e)=\sum_{\theta \subset e}\sum_{y\in h(\theta )\cap
D(\alpha )}\mathrm{I}(h(\theta ),L_{y})[y].
\end{equation*}%
Also, observe that
\begin{equation*}
\begin{array}{ll}
\lbrack \lbrack h(x_{11})]]=\epsilon ^{-a_{1}-a_{2}-a_{3}}[[y_{1}]], &
[[h(x_{12})]]=j\epsilon ^{-a_{1}}[[y_{1}]], \\
\lbrack \lbrack h(x_{21})]]=\epsilon ^{-a_{1}}[[y_{1}]], &
[[h(x_{22})]]=j\epsilon ^{-a_{1}-a_{2}-a_{3}}[[y_{1}]], \\
\lbrack \lbrack h(x_{31})]]=[[y_{2}]], & [[h(x_{32})]]=j\epsilon
^{-a_{1}-a_{2}}[[y_{2}]], \\
\lbrack \lbrack h(x_{41})]]=\epsilon ^{-a_{1}-a_{2}}[[y_{2}]], &
[[h(x_{42})]]=j[[y_{2}]].%
\end{array}%
\end{equation*}%
Now we are ready to calculate the obstruction cocycle $c_{Q_{4n}}^{\prime
}(h)(e)$.

\noindent (A) In the case $a_{1}=a_{2}=a_{3}=a_{4}$, according to previous
discussion, $c_{Q_{4n}}^{\prime }(h)(e)=[[h(x_{31})]]=[[y_{2}]].$

\noindent (B) The remaining cases of the computation of $[c_{Q_{4n}}^{\prime
}(h)(e)]$ are the following.

\noindent (B.1) If $a_{2}=a_{4}$ and $a_{1}\neq a_{3}$, then
\begin{equation*}
\lbrack c_{Q_{4n}}^{\prime
}(h)(e)]=((-1)^{(n+1)(a_{1}+a_{2}+a_{3})}(1+(-1)^{\binom{n}{2}%
})+(-1)^{(n+1)a_{1}}(1+(-1)^{\binom{n}{2}}))[[y_{1}]]+K.
\end{equation*}

\noindent (B.2) If $a_{2}=a_{4}\neq a_{1}=a_{3}$, then
\begin{equation*}
\lbrack c(h)(e)]=((-1)^{(n+1)(a_{1}+a_{2}+a_{3})}+(-1)^{(n+1)a_{1}+\binom{n%
}{2}})[[y_{1}]]+K.
\end{equation*}

\noindent (B.3) If $a_{1}=a_{3}$ and $a_{2}\neq a_{4}$, then
\begin{equation*}
\lbrack c_{Q_{4n}}^{\prime
}(h)(e)]=((-1)^{(n+1)(a_{1}+a_{2}+a_{3})}(1+(-1)^{\binom{n}{2}%
})+(-1)^{(n+1)a_{1}}(1+(-1)^{\binom{n}{2}}))[[y_{1}]]+K.
\end{equation*}%
In particular for $(a_{1},a_{2},a_{3},a_{4})=(p,p,p,2p)$%
\begin{equation*}
\lbrack c_{Q_{4n}}^{\prime
}(h)(e)]=2((1+(-1)^{\binom{5p}{2}}))[[y_{1}]]+K.
\end{equation*}

\noindent (B.4) If $a_{1}=a_{2}\neq a_{3}=a_{4}$, then
\begin{equation*}
\begin{array}{ll}
\lbrack c_{Q_{4n}}^{\prime }(h)(e)]= &
((-1)^{(n+1)(a_{1}+a_{2}+a_{3})}+(-1)^{(n+1)a_{1}+\binom{n}{2}})[[y_{1}]]-
\\
& ((-1)^{(n+1)(a_{1}+a_{2})}+1)(1+(-1)^{^{\binom{n}{2}}})[[y_{2}]]+K.%
\end{array}%
\end{equation*}

\noindent (B.5) If $a_{2}=a_{3}\neq a_{1}=a_{4}$, then
\begin{equation*}
\begin{array}{ll}
\lbrack c_{Q_{4n}}^{\prime }(h)(e)]= &
((-1)^{(n+1)(a_{1}+a_{2}+a_{3})}+(-1)^{(n+1)a_{1}})(1+(-1)^{^{\binom{n}{2}%
}})[[y_{1}]]- \\
& (1+(-1)^{(n+1)(a_{1}+a_{2})})[[y_{2}]]+K.%
\end{array}%
\end{equation*}

\noindent (B.6) In all the remaining cases we have
\begin{equation*}
\begin{array}{ll}
\lbrack c_{Q_{4n}}^{\prime }(h)(e)]= &
((-1)^{(n+1)(a_{1}+a_{2}+a_{3})}+(-1)^{(n+1)a_{1}})(1+(-1)^{^{\binom{n}{2}%
}})[[y_{1}]]- \\
& ((-1)^{(n+1)(a_{1}+a_{2})}+1)(1+(-1)^{^{\binom{n}{2}}})[[y_{2}]]+K.%
\end{array}%
\end{equation*}
\end{proof}

\subsubsection{The obstruction cocycle $c_{\mathbb{Z}_{n}}(h)$}

\noindent In this section we perform similar computation of the obstruction
cocycle for a $\mathbb{Z}_{n}$-map $h:S^{3}\rightarrow W_{n}$ in general
position, and complete the calculations originally started in \cite%
{VreZiv2002}. Here the arrangement $\mathcal{A}(\alpha )$ is the minimal $%
\mathbb{Z}_{n}$-invariant subspace arrangement generated by $L(\alpha )$
defined by (\ref{L(alfa)}).

\begin{theorem}
\label{Th.ObCycleZn}Let $c_{\mathbb{Z}_{n}}^{\prime }(h)\in C_{\mathbb{Z}%
_{n}}^{3}(S^{3},H_{2}(M;\mathbb{Z}))$ be the obstruction cocycle for the map
$h$ defined above. Let $K$ be the additive subgroup of $H_{2}(M;\mathbb{Z})$
generated by the elements of the form $g\cdot x-x$, $g\in \mathbb{Z}_{n}$, $%
x\in H_{2}(M;\mathbb{Z})$. Then the group of coinvariants is $H_{2}(M;%
\mathbb{Z})_{\mathbb{Z}_{n}}=H_{2}(M;\mathbb{Z})/K$.

(A) Let $(a_{1},a_{2},a_{3},a_{4})=(p,p,p,p)$. Then $[c_{\mathbb{Z}%
_{n}}^{\prime }(h)(d_{3})]=(-1)^{p}{[[}y_{1}]]+K={[[}y_{2}]]+K$.

(B) Let $(a_{1},a_{2},a_{3},a_{4})=(p,q,p,q)$. Then $[c_{\mathbb{Z}%
_{n}}^{\prime }(h)(d_{3})]=({(-1)}^{p}+(-1)^{q})[[y_{1}]]+K$.

(C) Let $(a_{1},a_{2},a_{3},a_{4})\notin \{(p,p,p,p),(p,q,p,q)\}$. Then
\begin{equation*}
\lbrack c_{\mathbb{Z}_{n}}^{\prime
}(h)(d_{3})]=((-1)^{(n+1)a_{1}}+(-1)^{(n+1)(a_{1}+a_{2}+a_{3})})[[y_{1}]]-((-1)^{(n+1)(a_{1}+a_{2})}+1)[[y_{2}]]+K
\end{equation*}

(D) Specially, let $(a_{1},a_{2},a_{3},a_{4})=(p,p,p,2p)$. Then
\begin{equation*}
\lbrack c_{\mathbb{Z}_{n}}^{\prime }(h)(d_{3})]=2[[y_{1}]]-2[[y_{2}]]+K.
\end{equation*}
\end{theorem}

\begin{proof}
The proof goes along the lines of the proof of theorems 4.1 and 4.2 in \cite%
{VreZiv2002}. We use again the formulas
\begin{equation*}
c_{\mathbb{Z}_{n}}(h)(\theta )=\sum_{y\in h(\theta )\cap D(\alpha )}\mathrm{I%
}(h(\theta ),L_{y})[y]\text{ and }c_{Q_{4n}}^{\prime
}(h)(d_{3})=\sum_{\theta \subset d_{3}}\sum_{y\in h(\theta )\cap D(\alpha )}%
\mathrm{I}(h(\theta ),L_{y})[y]\text{.}
\end{equation*}

\noindent (A) If $(a_{1},a_{2},a_{3},a_{4})=(p,p,p,p)$, then there is only
one simplex
\begin{equation*}
\epsilon ^{-a_{1}}\theta _{1}=\epsilon ^{-a_{1}-a_{2}}\theta _{2}=\epsilon
^{-a_{1}-a_{2}-a_{3}}\theta _{3}=\theta _{4}
\end{equation*}%
which belongs to the cell $d_{3}$ and intersects nontrivially the singular
set $h^{-1}(D(\alpha ))$. In this case there is only one element of the set $%
d_{3}\cap h^{-1}(D(\alpha ))$. Thus
\begin{equation*}
c_{\mathbb{Z}_{n}}^{\prime }(h)(d_{3})=\det (\epsilon ^{-a_{1}})\epsilon
^{-a_{1}}[[y_{1}]]=(-1)^{(4p+1)p}{[[}y_{1}]].
\end{equation*}%
\noindent (B) If $(a_{1},a_{2},a_{3},a_{4})=(p,q,p,q)$, there is only one
simplex which belongs to the cell $d_{3}$ and nontrivially intersects the
singular set $h^{-1}(\cup \mathcal{A}(\alpha ))$. Note that in this case $%
\left\vert d_{3}\cap h^{-1}(D(\alpha ))\right\vert =2$. This implies that
\begin{equation*}
c_{\mathbb{Z}_{n}}^{\prime }(h)(d_{3})=\det (\epsilon ^{-a_{1}})\epsilon
^{-a_{1}}[[y_{1}]]+\det (\epsilon ^{-a_{1}-a_{2}-a_{3}})\epsilon
^{-a_{1}-a_{2}-a_{3}}[[y_{3}]].
\end{equation*}%
\noindent (C) In all remaining cases $\left\vert d_{3}\cap h^{-1}(D(\alpha
))\right\vert =4$ and

$%
\begin{array}{ll}
c_{\mathbb{Z}_{n}}^{\prime }(h)(d_{3})= & \det (\epsilon ^{-a_{1}})\epsilon
^{-a_{1}}{[[}y_{1}]]-\det (\epsilon ^{-a_{1}-a_{2}})\epsilon ^{-a_{1}-a_{2}}{%
[[}y_{2}{]]}+ \\
& \det (\epsilon ^{-a_{1}-a_{2}-a_{3}})\epsilon
^{-a_{1}-a_{2}-a_{3}}[[y_{3}]]-{[[}y_{4}]].%
\end{array}%
$
\end{proof}

\subsection{Homology groups $H_{2}(M(\protect\alpha );\mathbb{Z})$ as $G$%
-modules}

\noindent Our objective in this section is to compute the torsion subgroups
of both $H_{2}(M;\mathbb{Z})_{\mathbb{Z}_{n}}$ and $H_{2}(M;\mathbb{Z}%
)_{Q_{4n}}$ (Proposition \ref{Th.ObCoCycleTorzioni}). In the first case $%
M=W_{n}-\cup \mathcal{A}(\alpha )$ is the complement of the minimal $\mathbb{%
Z}_{n}$-arrangement $\mathcal{A}(\alpha )$ containing subspace $L(\alpha )$,
and in the second case $\mathcal{A}(\alpha )$ is the minimal $D_{2n}$%
-arrangement containing $L(\alpha )$ (defined by (\ref{L(alfa)})).

\noindent The idea is to use Poincar\'{e}-Alexander duality and work with
the arrangement $\mathcal{A}(\alpha )$ instead of the complement $M$. The
complement $M=W_{n}-\cup \mathcal{A}(\alpha )$ can be written in the form $%
S^{n-1}-\cup \widehat{\mathcal{A}}(\alpha )$, where $\widehat{\mathcal{A}}%
(\alpha )$ is the compactification of the arrangement $\mathcal{A}(\alpha )$%
. There is a sequence of isomorphisms
\begin{eqnarray*}
H_{2}(M,\mathbb{Z}) &=&H_{2}(S^{n-1}-\cup \widehat{\mathcal{A}}(\alpha ),%
\mathbb{Z})\text{ \ \ \ Poincar\'{e}-Alexander duality} \\
&\cong &H^{(n-1)-2-1}(\cup \widehat{\mathcal{A}}(\alpha ),\mathbb{Z})\text{
\ \ \ Universal Coefficient theorem} \\
&\cong &\mathrm{Hom}(H_{n-4}(\cup \widehat{\mathcal{A}}(\alpha ),\mathbb{Z}),%
\mathbb{Z})\oplus \mathrm{Ext}(H_{n-3}(\cup \widehat{\mathcal{A}}(\alpha ),%
\mathbb{Z}),\mathbb{Z}).
\end{eqnarray*}%
The isomorphism of universal coefficient theorem is a $\mathbb{Z}_{n}$,
respectively a $Q_{4n}$-map,but the Poincar\'{e}-Alexander duality map is a $%
\mathbb{Z}_{n}$, i.e. $Q_{4n}$-map up to a orientation character. Let $o$ be
an orientation of the $\mathbb{Z}_{n}$, respectively $Q_{4n}$-sphere $%
S^{n-1} $. Then $o$ determines Poincar\'{e}-Alexander duality map $\gamma
_{o}$ \cite{Mu}. In particular, the mapping is $Q_{4n}$-equivariant up to
the orientation $g\cdot o=\det (g)\cdot o$, where $g\in \mathbb{Z}%
_{n}\subseteq GL_{n}(\mathbb{R})$, i.e. $g\in Q_{4n}\subseteq GL_{n}(\mathbb{%
R})$. Since the maximal elements of the arrangement $\mathcal{A}(\alpha )$
are $(n-4)$-dimensional linear subspaces, it follows that $H_{n-3}(\cup
\widehat{\mathcal{A}}(\alpha ),\mathbb{Z})=0$ and $\mathrm{Ext}(H_{n-3}(\cup
\widehat{\mathcal{A}}(\alpha ),\mathbb{Z}),\mathbb{Z})=0$. Hence,
\begin{equation}
H_{2}(M(\mathcal{A}(\alpha )),\mathbb{Z})\cong \mathrm{Hom}(H_{n-4}(\cup
\widehat{\mathcal{A}}(\alpha ),\mathbb{Z}),\mathbb{Z}).
\label{EqIzomorfizam-Hom}
\end{equation}%
The Ziegler-\v{Z}ivaljevi\'{c} formula implies the following decomposition
(assuming $\mathbb{Z}$ as coefficients)
\begin{eqnarray*}
H_{n-4}(\cup \widehat{\mathcal{A}}(\alpha )) &\cong &\underset{p\in P(\alpha
)}{\bigoplus }H_{n-4}(\Delta (P(\alpha )_{<p})\ast S^{\dim p})\cong \underset%
{p\in P(\alpha )}{\bigoplus }H_{n-4}(\Sigma ^{\dim p+1}(\Delta (P(\alpha
)_{<p}))\mathbf{)} \\
&\cong &\underset{p\in P(\alpha )}{\bigoplus }\tilde{H}_{n-5-\dim p}(\Delta
(P(\alpha )_{<p})\mathbf{)}\cong \underset{d=0}{\overset{n-4}{\bigoplus }}%
\underset{p\in P(\alpha ):\dim p=d}{\bigoplus }\tilde{H}_{n-5-d}(\Delta
(P(\alpha )_{<p})\mathbf{).}
\end{eqnarray*}%
where $P(\alpha )$ is the intersection poset of the arrangement $\mathcal{A}%
(\alpha )$. Thus, in both cases $\mathbb{Z}_{n}$ and $Q_{4n}$ we have to
determine $(\forall p\in P(\alpha ))\,\tilde{H}_{n-5-\dim p}(\Delta
(P(\alpha )_{<p}))$ with $\tilde{H}_{-1}(\emptyset )=\mathbb{Z}$. Observe
that the dimension of the simplicial complex $\Delta (P(\alpha )_{<p})$ is
less or equal then $n-5-\dim p$ and $\tilde{H}_{n-5-\dim p}(\Delta (P(\alpha
))_{<p})\neq 0$ implies that there must be at least one chain of length $\
n-5-\dim p.$

\subsubsection{Case 1: $\mathcal{A}(\protect\alpha )$ as the minimal $%
\mathbb{Z}_{n}$-invariant subspace arrangement}

\begin{theorem}
\label{Th.HomologijaZn}Let $M=W_{n}-\cup \mathcal{A}(\alpha )$ be the
complement of the minimal $\mathbb{Z}_{n}$-arrangement $\mathcal{A}(\alpha )$
containing subspace $L(\alpha )$ defined by \ref{L(alfa)}. Let $\alpha
=(a_{1},a_{2},a_{3},a_{4})\in \mathbb{N}^{4}$ and $a_{1}+a_{2}+a_{3}+a_{4}=n$%
. Then

\noindent (A) If $\alpha =(p,p,p,p)$, then $H_{2}(M;\mathbb{Z})\cong \mathbb{%
Z[Z}_{n}/\epsilon ^{p}\mathbb{Z}_{n}\mathbb{]}\cong \mathbb{Z}^{n/4}$;

\noindent (B) If $\alpha =(p,q,p,q)$ and $p\neq q$, then $H_{2}(M;\mathbb{Z}%
)\cong \mathbb{Z[Z}_{n}/\epsilon ^{p+q}\mathbb{Z}_{n}\mathbb{]}\cong \mathbb{%
Z}^{n/2}$;

\noindent (C) If $\alpha =(p,p,p,2p)$, then $H_{2}(M;\mathbb{Z})\cong
\mathbb{Z[Z}_{n}\mathbb{]\oplus }\frac{n}{5}\mathbb{Z}^{4}\cong \mathbb{Z}%
^{n+\frac{4n}{5}}$;

\noindent (D) If $\alpha =(1,1,1,3)$, then $H_{2}(M;\mathbb{Z})\cong \mathbb{%
Z[Z}_{6}\mathbb{]\oplus Z[Z}_{6}\mathbb{]\oplus Z}\cong \mathbb{Z}^{13}$;

\noindent (E) If $\alpha =(q,q,q,p)$, $p\neq q$, $n>6$, then $H_{2}(M;%
\mathbb{Z})\cong \mathbb{Z[Z}_{n}\mathbb{]\oplus Z[Z}_{n}\mathbb{]}\cong
\mathbb{Z}^{2n}$;

\noindent (F) If $\alpha =(p,q,p,p+q)$, $p\neq q$, then $H_{2}(M;\mathbb{Z}%
)\cong \mathbb{Z[Z}_{n}\mathbb{]\oplus Z[Z}_{n}\mathbb{]}\cong \mathbb{Z}%
^{2n}$;

\noindent (G) If $\alpha =(p,q,p+q,p+q)$, then $H_{2}(M;\mathbb{Z})\cong
\mathbb{Z[Z}_{n}\mathbb{]\oplus Z[Z}_{n}\mathbb{]\oplus Z[Z}_{n}/\epsilon
^{p+q}\mathbb{Z}_{n}\mathbb{]}$;

\noindent (H) In all the remaining cases $H_{2}(M;\mathbb{Z})\cong \mathbb{%
Z[Z}_{n}\mathbb{]}$.
\end{theorem}

\noindent Before discussing details of the proof let us make a few general
observations. The computation of $H_{2}(M;\mathbb{Z})$ will rely on the
isomorphism (\ref{EqIzomorfizam-Hom}) and the fact (which will be proved)
that the homology $H_{n-4}(\cup \widehat{\mathcal{A}}(\alpha ),\mathbb{Z})$
is free. The computation of the group $H_{n-4}(\cup \widehat{\mathcal{A}}%
(\alpha ),\mathbb{Z})$ is based on the Z-\v{Z} formula and again we will
have to discuss each of the cases separately. The Z-\v{Z} formula implies
that we actually compute $\tilde{H}_{\dim L-\dim p-1}(\Delta (P(\alpha
))_{<p};\mathbb{Z})$ for each $p\in P(\alpha )$.

\begin{proof}
\textbf{(A)} In this case there are $n/4$ maximal elements $L$, $\epsilon L$%
, ..., $\epsilon ^{p-1}L$ of the arrangement $\mathcal{A}(\alpha )$. Since
for every $i,j\in \{0,..,p-1\}$, $i\neq j$, \ \ \
\begin{equation*}
\dim L-\dim (\epsilon ^{i}L\cap \epsilon ^{j}L)>1
\end{equation*}
the Z-\v{Z} formula implies that only maximal elements contribute to $%
H_{n-4}(\cup \widehat{\mathcal{A}}(\alpha ),\mathbb{Z})$.

\textbf{(B) }and \textbf{(H) }The argument for these cases is exactly the
same as for (A), except for the number of maximal elements of the
arrangement $\mathcal{A}(\alpha )$.

\textbf{(C) }There are four maximal elements $\epsilon ^{p}L$, $\epsilon
^{2p}L$, $\epsilon ^{3p}L$ and $\epsilon ^{4p}L$ with the $\dim (L\cap
\epsilon ^{k\cdot p}L)=\dim L-1$ for $k=1,..,4$. In addition
\begin{equation*}
L\cap \epsilon ^{p}L=L\cap \epsilon ^{2p}L=L\cap \epsilon ^{3p}L=L\cap
\epsilon ^{4p}L\text{.}
\end{equation*}%
Let $V=L\cap \epsilon ^{q}L\cap \epsilon ^{2q}L\cap \epsilon
^{3q}L\cap \epsilon ^{4q}L$. Then the Hasse diagram of the
subposet $\{p\in P(\alpha ):\dim p\geq n-5\}$ is depicted in the
Figure~\ref{fig:sl4}.

\begin{figure}[htb]
\centering
\includegraphics[scale=0.80]{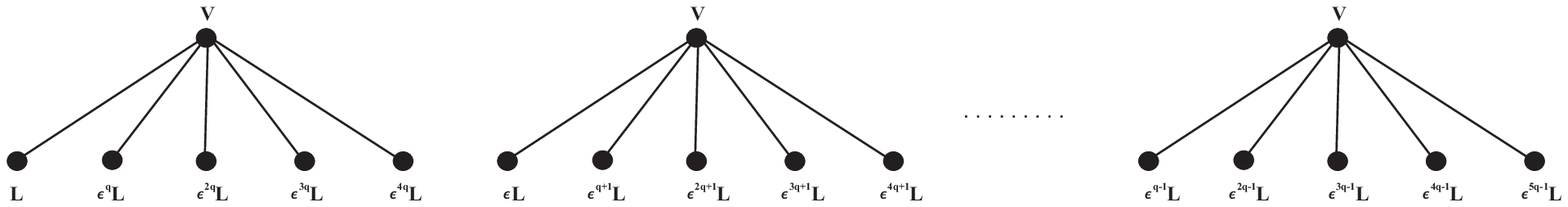}
\caption{The Hasse diagram for the case (C)} \label{fig:sl4}
\end{figure}

\noindent Before using the Z-\v{Z} decomposition observe that for $n\geq 5$%
\begin{equation*}
(\forall p\in P(\alpha ))\,\dim p\leq n-6\,\Longrightarrow \,\tilde{H}%
_{n-5-\dim p}(\Delta (P(\alpha )_{<p}))=0.
\end{equation*}

\noindent Indeed, for each element $q\in P(\alpha )_{<p}$ such that $\dim
q=n-4$ there exists a unique element $p_{q}\in P(\alpha )_{<p}$ with the
property $\dim p_{q}=n-5$ and $q<p_{q}$. There is a monotone map $f:P(\alpha
)_{<p}\rightarrow P(\alpha )_{<p}-\{q|\dim q=n-4\}$ defined by
\begin{equation*}
\begin{array}{ccc}
q & \longmapsto & \left\{
\begin{array}{cc}
q, & \text{for }\dim q\leq n-5 \\
p_{q} & \text{for }\dim q=n-4%
\end{array}%
\right.%
\end{array}%
\end{equation*}%
which satisfies conditions of the Quillen fiber lemma. This implies that $f$
induces a homotopy equivalence, hence
\begin{equation*}
\tilde{H}_{n-5-\dim p}(\Delta (P(\alpha )_{<p}))=\tilde{H}_{n-5-\dim
p}(\Delta (Q))=0
\end{equation*}%
since $\dim \Delta (Q)<n-5-\dim p$.

\noindent Thus, the only relevant part of the poset $P(\alpha )$ for the
computation of $H_{n-4}(\cup \widehat{\mathcal{A}}(\alpha ),\mathbb{Z})$ is
the part in the above picture.

\textbf{(D)} The proof follows from the Z-\v{Z} decomposition and
the Hasse diagram of the intersection poset, shown in the
Figure~\ref{fig:sl2}.

\begin{figure}[htb]
\centering
\includegraphics[scale=0.80]{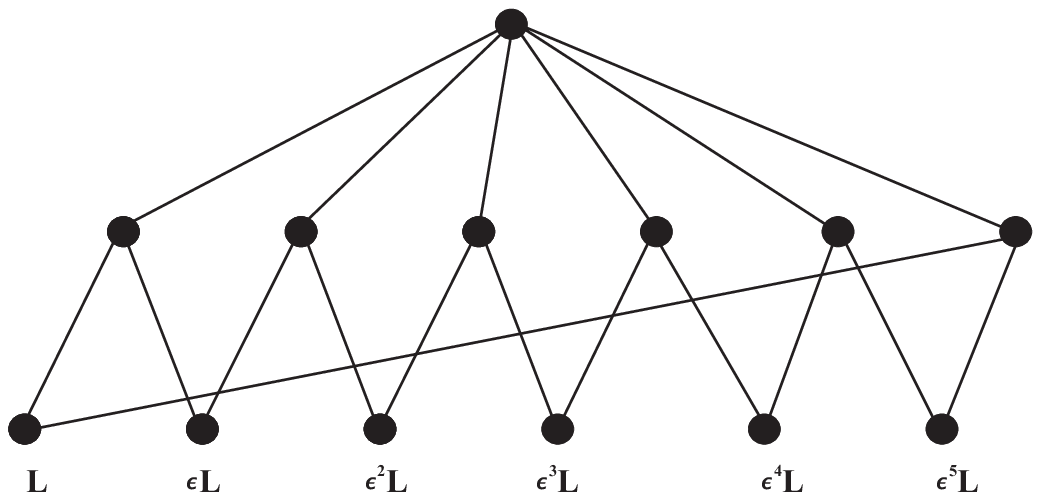}
\caption{The Hasse diagram for the case (D)} \label{fig:sl2}
\end{figure}

\textbf{(E)} When $\alpha =(a_{1},a_{2},a_{3},a_{4})=(q,q,q,p)$, $p\neq 2q$,
there are exactly two maximal elements $\epsilon ^{q}L$ and $\epsilon
^{2q+p}L$ in $P(\alpha )$ with the maximal $L$ intersection, $\dim (L\cap
\epsilon ^{q}L)=\dim (L\cap \epsilon ^{2q+p}L)=\dim L-1$. It is not hard to
describe the first three levels of the intersection poset $P(\alpha )$.
There are $n$ elements $L,\epsilon L,..,\epsilon ^{n-1}L$ of dimension $n-4$%
, $n$ elements
\begin{eqnarray*}
&&L\cap \epsilon ^{q}L,\epsilon ^{q}L\cap \epsilon ^{2q}L,...,\epsilon
^{2q+p}L\cap L;\,\epsilon L\cap \epsilon ^{q+1}L,\epsilon ^{q+1}L\cap
\epsilon ^{2q+1}L,...,\epsilon ^{2q+p+1}L\cap \epsilon L; \\
&&\epsilon ^{d-1}L\cap \epsilon ^{q+d-1}L,\epsilon ^{q+d-1}L\cap \epsilon
^{2q+d-1}L,...,\epsilon ^{2q+p+d-1}L\cap \epsilon ^{d-1}L
\end{eqnarray*}%
of dimension $n-5$ (where $d=(q,n)$) and finally $n$ elements
\begin{equation*}
L\cap \epsilon ^{q}L\cap \epsilon ^{2q+p}L,\epsilon (L\cap \epsilon
^{q}L\cap \epsilon ^{2q+p}L),..,\epsilon ^{n-1}(L\cap \epsilon ^{q}L\cap
\epsilon ^{2q+p}L)
\end{equation*}%
of dimension $n-6$. Thus, the Hasse diagram of the subposet
$\{p\in P(\alpha ):\dim p\geq n-5\}$ is of the form depicted in
the Figure~\ref{fig:sl1}.

\begin{figure}[htb]
\centering
\includegraphics[scale=1.00]{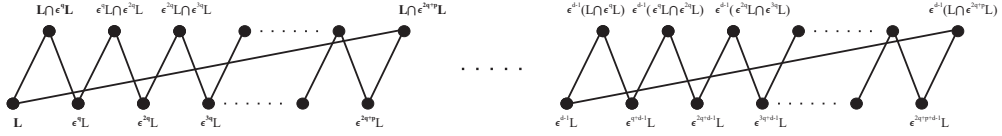}
\caption{The first Hasse diagram for the case (E)} \label{fig:sl1}
\end{figure}

\noindent We prove that only this part of the Hasse diagram of the
intersection poset $P(\alpha )$ contributes to the homology $H_{n-4}(\cup
\widehat{\mathcal{A}}(\alpha ),\mathbb{Z})$. Indeed,\textit{\ }for $n>6$\
\begin{equation*}
(\forall p\in P(\alpha ))\,\dim p\leq n-6\,\Longrightarrow \,\tilde{H}%
_{n-5-\dim p}(\Delta (P(\alpha )_{<p}))=0.
\end{equation*}%
\noindent Observe that $\dim \Delta (P(\alpha )_{<p})\leq n-5-\dim p$. If $%
\dim \Delta (P(\alpha )_{<p})<n-5-\dim p$ we have nothing to do, so we
assume that $\dim \Delta (P(\alpha )_{<p})=n-5-\dim p$. Let $K=\{q\in
P(\alpha )_{<p}|\dim q\neq n-6\}$. Let us show that the inclusion $%
i:K\hookrightarrow P(\alpha )_{<p}$ satisfies the assumptions of
the Quillen fiber lemma. It suffices to check that for every
$q=\epsilon ^{i}(L\cap \epsilon ^{q}L\cap \epsilon ^{2q+p}L)\in
P(\alpha )_{<p}-K$ the order complex $\Delta (i^{-1}((P(\alpha
)_{<p})_{\leq q}))$ is contractible. For such a $q=\epsilon
^{i}(L\cap \epsilon ^{q}L\cap \epsilon ^{2q+p}L)$, the Hasse
diagram of the poset $i^{-1}((P(\alpha )_{<p})_{\leq q})$ is
depicted in the Figure~\ref{fig:sl3} we observe that $\Delta
(i^{-1}((P(\alpha )_{<p})_{\leq q}))$ is obviously contractible.
Hence, $\Delta (P(\alpha
)_{<p})\simeq \Delta (K)$. Since, $\dim \Delta (P(\alpha )_{<p})<n-5-\dim p$%
, we conclude that $\tilde{H}_{n-5-\dim p}(K)=0$.

\begin{figure}[htb]
\centering
\includegraphics[scale=0.70]{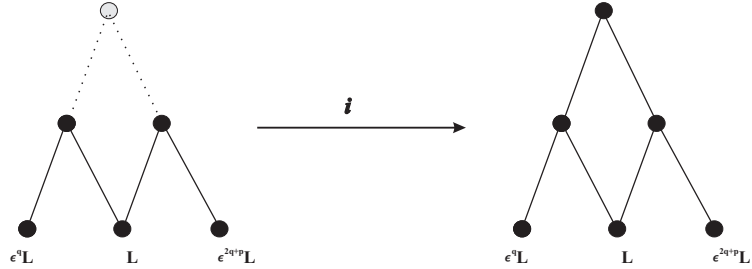}
\caption{The second Hasse diagram for the case (E)}
\label{fig:sl3}
\end{figure}

\noindent As before, from the Ziegler-\v{Z}ivaljevi\'{c} formula we deduce
the following decomposition
\begin{equation*}
H_{n-4}(\cup \widehat{\mathcal{A}}(\alpha );\mathbb{Z})\cong \underset{d=n-5}%
{\overset{n-4}{\bigoplus }}\underset{p\in P(\alpha ):\dim p=d}{\bigoplus }%
\tilde{H}_{n-5-d}(\Delta (P(\alpha )_{<p});\mathbb{Z}\mathbf{)}\cong \mathbb{%
Z}[\mathbb{Z}_{n}]\oplus \mathbb{Z}[\mathbb{Z}_{n}]\mathbf{.}
\end{equation*}

\textbf{(F)} There are exactly two maximal elements $\epsilon ^{p+q}L$ and $%
\epsilon ^{2p+q}L$ with the property
\begin{equation*}
\dim (L\cap \epsilon ^{p+q}L)=\dim (L\cap \epsilon ^{2p+q}L)=\dim L-1\text{.}
\end{equation*}%
Since the subposet $\{p\in P(\alpha ):\dim p\geq n-5\}$ of $P(\alpha )$ has
the same shape as the appropriate subposet in the proof of the part (E), it
is not hard to see that the computation for this case goes completely along
the lines of the preceding proof.

\textbf{(G) }In this case there are two maximal elements $\epsilon ^{p+q}L$
and $\epsilon ^{2p+2q}L$ with the maximal $L$ intersection , i.e.
\begin{equation*}
\dim (L\cap \epsilon ^{p+q}L)=\dim (L\cap \epsilon ^{2p+2q}L)=\dim L-1\text{.%
}
\end{equation*}%
The elements of $P(\alpha )$ of dimension $n-5$ are: $L\cap
\epsilon ^{p+q}L,\epsilon (L\cap \epsilon ^{p+q}L),..,\epsilon
^{n-1}(L\cap \epsilon ^{p+q}L)$. There are $p+q$ elements of
dimension $n-6$ for $n>6$: $L\cap \epsilon ^{p+q}L\cap \epsilon
^{2p+2q}L,\,\epsilon (L\cap \epsilon ^{p+q}L\cap \epsilon
^{2p+2q}L),..,\,\epsilon ^{p+q-1}(L\cap \epsilon ^{p+q}L\cap
\epsilon ^{2p+2q}L)$ and only one for $n=6$. We use this
information to draw the Hasse diagram (Figure~\ref{fig:sl5}) for
the subposet $\{p\in P(\alpha ):\dim p\geq n-6\}$.

\begin{figure}[htb]
\centering
\includegraphics[scale=0.70]{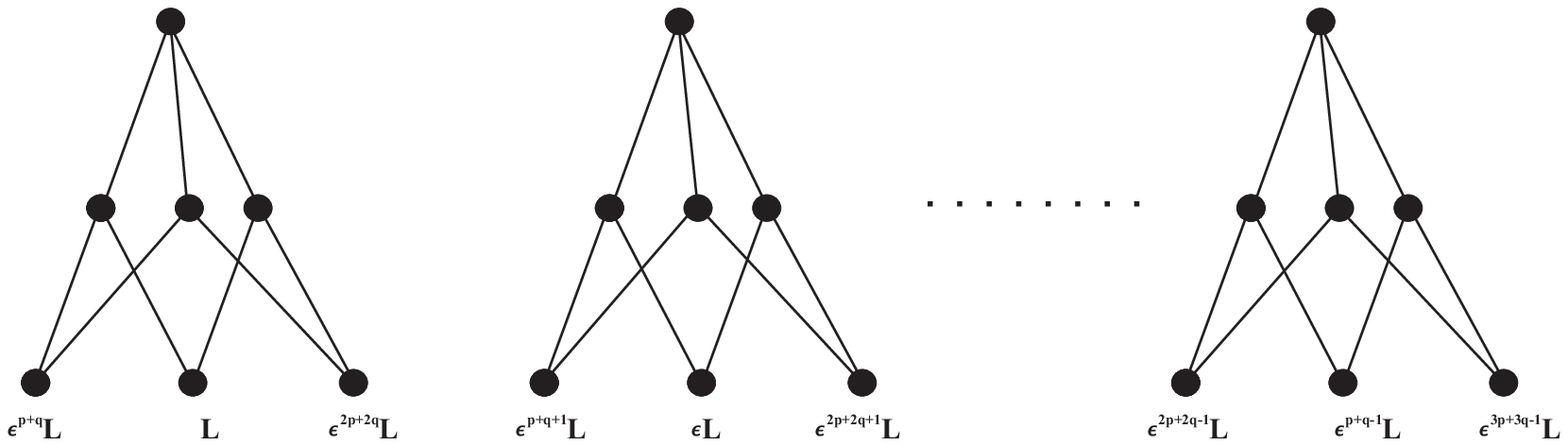}
\caption{The Hasse diagram for the case (G)} \label{fig:sl5}
\end{figure}

Let us convince ourselves that for $n>6$
\begin{equation*}
(\forall p\in P(\alpha ))\,\dim p\leq n-7\Longrightarrow \tilde{H}_{n-5-\dim
p}(\Delta (P(\alpha )_{<p}))=0.
\end{equation*}

\noindent\ Let $R\overset{def}{=}P(\alpha )_{<p}\cap \{p\in P(\alpha )|\dim
p=n-6\}=\{p_{1},..,p_{k}\}$ where $k\leq p+q$ and $Q=P(\alpha )_{<p}\cap
\{q\in P(\alpha )|\dim q\leq n-6\}$. Since for every element $q\in P(\alpha
)_{<p}\cap \{p\in P(\alpha )|\dim p>n-6\}$ there exists an unique element $%
p_{q}\in R$ such that $q<p_{q},$ we can define a monotone map $f:P(\alpha
)_{<p}\rightarrow Q$ by the formula
\begin{equation*}
f(q)=\left\{
\begin{array}{cc}
q & \text{, }\dim q\leq n-6 \\
p_{q} & \text{, }\dim q>n-6%
\end{array}%
\right.
\end{equation*}%
Again the Quillen fiber lemma implies that $f$ induces a homotopy
equivalence $\Delta (P(\alpha )_{<p})\simeq \Delta (Q)$. Indeed, for each $%
q\in Q$ the poset $f^{-1}((P(\alpha )_{<p})_{\leq q})$ has a maximum, hence $%
\Delta (f^{-1}((P(\alpha )_{<p})_{\leq q}))$ is contractible. It follows
from $\dim \Delta (Q)<n-5-\dim p$ that
\begin{equation*}
\tilde{H}_{n-5-\dim p}(\Delta (P(\alpha )_{<p}))=\tilde{H}_{n-5-\dim
p}(\Delta (Q))=0.
\end{equation*}

\noindent From \ here we infer the that Ziegler-\v{Z}ivaljevi\'{c} formula
decomposition reduces to the first three levels of the Hasse diagram of $%
P(\alpha )\ $and
\begin{equation*}
H_{n-4}(\cup \widehat{\mathcal{A}}(\alpha ),\mathbb{Z})\cong \mathbb{Z[Z}_{n}%
\mathbb{]\oplus Z[Z}_{n}\mathbb{]\oplus Z[Z}_{n}/\epsilon ^{p+q}\mathbb{Z}%
_{n}\mathbb{]}.
\end{equation*}

\textbf{(H) }In all the remaining cases there are exactly $n$ maximal
elements of the arrangement, an no intersection of two maximal elements is
of dimension $n-5$. Therefore, only maximal elements contribute to the
homology $H_{n-4}(\cup \widehat{\mathcal{A}}(\alpha ),\mathbb{Z})$.
\end{proof}

\subsubsection{Case 2: $\mathcal{A}(\protect\alpha )$ as the minimal $Q_{4n}$%
-invariant subspace arrangement}

\begin{theorem}
\label{Th.HomologijaQ4n}Let $M=W_{n}-\cup \mathcal{A}(\alpha )$ be the
complement of the minimal $D_{2n}$-invariant arrangement $\mathcal{A}(\alpha
)$ containing the subspace $L(\alpha )$ defined by (\ref{L(alfa)}) (or
equivalently the $Q_{4n}$-arrangement, where $j^{2}$ acts trivially). Let $%
\alpha =(a_{1},a_{2},a_{3},a_{4})\in \mathbb{N}^{4}$ and $%
a_{1}+a_{2}+a_{3}+a_{4}=n$. Then

\noindent (A) If $\alpha =(p,p,p,p)$, then $H_{2}(M;\mathbb{Z})\cong \mathbb{%
Z[Z}_{n}/\epsilon ^{p}\mathbb{Z}_{n}\mathbb{]}\cong \mathbb{Z}^{n/4}$;

\noindent (B) If $\alpha =(p,q,p,q)$ and $p\neq q$, then $H_{2}(M;\mathbb{Z}%
)\cong \mathbb{Z[Z}_{n}/\epsilon ^{p+q}\mathbb{Z}_{n}\mathbb{]}\cong \mathbb{%
Z}^{n/2}$;

\noindent (C) If $\alpha =(p,p,p,2p)$, then $H_{2}(M;\mathbb{Z})\cong
\mathbb{Z[Z}_{n}\mathbb{]\oplus }\frac{n}{5}\mathbb{Z}^{4}\cong \mathbb{Z}%
^{n+\frac{4n}{5}}$;

\noindent (D) If $\alpha =(1,1,1,3)$, then $H_{2}(M;\mathbb{Z})\cong \mathbb{%
Z[Z}_{6}\mathbb{]\oplus Z[Z}_{6}\mathbb{]\oplus Z}\cong \mathbb{Z}^{13}$;

\noindent (E) If $\alpha =(q,q,q,p)$ and $p\notin \{q,2q\}$, then $H_{2}(M;%
\mathbb{Z})\cong \mathbb{Z[Z}_{n}\mathbb{]\oplus Z[Z}_{n}\mathbb{]}\cong
\mathbb{Z}^{2n}$;

\noindent (F) If $\alpha =(p,q,p,p+q)$ and $p\neq q$, then $H_{2}(M;\mathbb{Z%
})\cong \mathbb{Z[Z}_{n}\mathbb{]\oplus Z[Z}_{n}\mathbb{]}\cong \mathbb{Z}%
^{2n}$;

\noindent (G) If $\alpha =(p,p,2p,2p)$, then $H_{2}(M;\mathbb{Z})\cong
\mathbb{Z[Z}_{n}\mathbb{]\oplus Z[Z}_{n}\mathbb{]\oplus Z[Z}_{n}/\epsilon
^{2p}\mathbb{Z}_{n}\mathbb{]}$;

\noindent (H) If $\alpha =(p,q,p+q,p+q)$ and $p\neq q$, then $H_{2}(M;%
\mathbb{Z})\cong ?$.

\noindent (I) If $\alpha =(p,p,q,q)$ and $p\neq q$, $p\neq 2q$, $q\neq 2p$,
then $H_{2}(M;\mathbb{Z})\cong \mathbb{Z[Z}_{n}\mathbb{]}$;

\noindent (J) If $\alpha =(p,q,r,q)$ and $r\notin \{q,p,p+q\}$, then $%
H_{2}(M;\mathbb{Z})\cong \mathbb{Z[Z}_{n}\mathbb{]}$;

\noindent (K) If $\alpha =(p,q,r,p+q)$ and $r\notin \{p,p+q\}$, then $%
H_{2}(M;\mathbb{Z})\cong ?$;

\noindent (L) If $\alpha =(p,p,q,q+r)$, and $p\neq q$, $r\neq 0$, then $%
H_{2}(M;\mathbb{Z})\cong ?$;

\noindent (M) In all the remaining cases $H_{2}(M;\mathbb{Z})\cong \mathbb{Z[%
}D_{2n}\mathbb{]}$.
\end{theorem}

\begin{remark}
The computation of $H_{2}(M;\mathbb{Z})$ in the cases (H), (K) and (L) is
probably within reach of the existing methods, however we leave it as an
open problem.
\end{remark}

\noindent Like in the preceding section the computation of $H_{2}(M;\mathbb{Z%
})$ will rely on the isomorphism (\ref{EqIzomorfizam-Hom}) and the fact that
the homology $H_{n-4}(\cup \widehat{\mathcal{A}}(\alpha ),\mathbb{Z})$ is
free. This means that we actually compute $H_{n-4}(\cup \widehat{\mathcal{A}}%
(\alpha ),\mathbb{Z})$, i.e. with the help of the Z-\v{Z} formula we compute
$\tilde{H}_{\dim L-\dim p-1}(\Delta (P(\alpha ))_{<p};\mathbb{Z})$ for each $%
p\in P(\alpha )$.

\begin{proof}
\textbf{(A)-(G) and (I)-(J) }In these cases minimal $D_{2n}$-arrangement
containing the subspace $L(\alpha )$ coincides with the minimal $\mathbb{Z}%
_{n}$-arrangement containing the same subspace $L(\alpha )$. Thus the
assertions are consequences of Theorem \ref{Th.HomologijaZn}.

\textbf{(M) }In all the remaining cases there are exactly $2n$ maximal
elements of the arrangement, an no intersection of two maximal elements is
of dimension $n-5$. Therefore, only maximal elements contribute to the
homology $H_{n-4}(\cup \widehat{\mathcal{A}}(\alpha ),\mathbb{Z})$.
\end{proof}

\subsection{The coinvariants of $H_{2}(M(\protect\alpha );\mathbb{Z})$}

\noindent Relying on the results from the previous section we are able to
read off the $\mathbb{Z}_{n}$, respectively $Q_{4n},$ coinvariants of the
module $H_{2}(M;\mathbb{Z})$.

\subsubsection{Case 1: $\mathbb{Z}_{n}$-coinvariants of $H_{2}(M;\mathbb{Z})$%
}

\begin{theorem}
\label{Th.KoinvarijanteZn}Let $M=W_{n}-\cup \mathcal{A}(\alpha )$ be the
complement of the minimal $\mathbb{Z}_{n}$-arrangement $\mathcal{A}(\alpha )$
containing the subspace $L(\alpha )$ defined by \ref{L(alfa)}. Let $\alpha
=(a_{1},a_{2},a_{3},a_{4})\in \mathbb{N}^{4}$ and $a_{1}+a_{2}+a_{3}+a_{4}=n$%
. Then%
\begin{equation*}
\begin{tabular}{|l|l|l|}
\hline
& $\alpha $ & $H_{2}(M;\mathbb{Z})_{\mathbb{Z}_{n}}$ \\ \hline
(A) & $(p,p,p,p)$ & $\mathbb{Z}_{2}$ \\ \hline
(B) & $(p,q,p,q)$, $p\neq q$ & $\mathbb{Z}$ \\ \hline
(C) & $(p,p,p,2p)$ & $\mathbb{Z\oplus Z}_{5}$ \\ \hline
(D) & $(1,1,1,3)$ & $\mathbb{Z\oplus Z\oplus Z}_{2}$ \\ \hline
(E) & $(q,q,q,p)$, $p\notin \{q,2q\}$ & $\mathbb{Z\oplus Z}$ \\ \hline
(F) & $(p,q,p,p+q)$, $p\neq q$ & $\mathbb{Z\oplus Z}$ \\ \hline
(G) & $(p,q,p+q,p+q)$ & $\mathbb{Z\oplus Z\oplus Z}$ \\ \hline
(H) & In all the remaining cases & $\mathbb{Z}$ \\ \hline
\end{tabular}%
\end{equation*}
\end{theorem}

\noindent To prove this theorem we have to keep in mind that the of Poincar%
\'{e}-Alexander duality map is an equivariant up to an orientation
character. To compute $\mathbb{Z}_{n}$-coinvariants $H_{2}(M(\mathcal{A}%
(\alpha )),\mathbb{Z})_{\mathbb{Z}_{n}}$\ we pass to the $\mathbb{Z}_{n}$%
-module $H_{n-4}(\cup \widehat{\mathcal{A}}(\alpha ),\mathbb{Z})$, but with
the modified action. More precisely, let $l\in H_{n-4}(\cup \widehat{%
\mathcal{A}}(\alpha ),\mathbb{Z})$ and $g\in \mathbb{Z}_{n}$, then
\begin{equation*}
g\cdot _{\text{m}}l=\det (g)~g\cdot l
\end{equation*}%
where $\cdot _{\text{m}}$ is the new modified action, and $\cdot $ the old
action. Let us denote by $\sim $ the equivalence relation on $H_{n-4}(\cup
\widehat{\mathcal{A}}(\alpha ),\mathbb{Z})$ such that the subgroup of
elements equivalent to zero is generated by the elements of the form $g\cdot
_{\text{m}}x-x$, $g\in \mathbb{Z}_{n}$, $x\in H_{n-4}(\cup \widehat{\mathcal{%
A}}(\alpha ),\mathbb{Z})$. In other words there is an isomorphism $H_{2}(M;%
\mathbb{Z})_{\mathbb{Z}_{n}}\cong H_{n-4}(\cup \widehat{\mathcal{A}}(\alpha
),\mathbb{Z})/\sim $.

\begin{proof}
\textbf{(A)} Let $\alpha =(p,p,p,p)$ and $l\in H_{n-4}(\cup \widehat{%
\mathcal{A}}(\alpha ),\mathbb{Z})$ be the homology class induced by the
subspaces $L$. The equality $L=\epsilon ^{p}\cdot L$ translates to the
equality $l=(-1)^{p-1}\epsilon ^{p}\cdot l$ in homology, because $\det
(\epsilon ^{p}|_{L})=(-1)^{p-1}$. Now in coinvariants we have
\begin{equation*}
l\sim \epsilon ^{p}\cdot _{\text{m}}l=\det (\epsilon ^{p})~\epsilon
^{p}\cdot l=(-1)^{p(n+1)}~\epsilon ^{p}\cdot l=(-1)^{p(4p+1)}(-1)^{p-1}l=-l
\end{equation*}%
which implies the first result.

\noindent \textbf{(B) }For $\alpha =(p,q,p,q)$ let $l$ be the homology class
induced by the subspaces $L$. The equality $L=\epsilon ^{p+q}\cdot L$
becomes relation $l=(-1)^{p+q}\epsilon ^{p+q}\cdot l$ in homology ($\det
(\epsilon ^{p}|_{L})=(-1)^{p+q}$), and in coinvariants%
\begin{equation*}
l\sim \epsilon ^{p+q}\cdot _{\text{m}}l=\det (\epsilon ^{p+q})~\epsilon
^{p+q}\cdot l=(-1)^{(p+q)(n+1)}~\epsilon ^{p+q}\cdot
l=(-1)^{(p+q)(2p+2q+1)}(-1)^{p+q}l=l\text{.}
\end{equation*}%
\qquad \noindent\

\noindent \textbf{(C) }Let $l,k\in H_{n-4}(\cup \widehat{\mathcal{A}}(\alpha
),\mathbb{Z})$ be homology classes induced by subspaces $L$ and $L\cap
\epsilon ^{p}L$. Then $l$ generates $\mathbb{Z[Z}_{n}\mathbb{]}$ part in $%
H_{n-4}(\cup \widehat{\mathcal{A}}(\alpha ),\mathbb{Z})$ and there are no
nontrivial relations on $l$. For $k$ there is only one relevant relation.
The homology equality $k+\epsilon ^{p}\cdot k+\epsilon ^{2p}\cdot k+\epsilon
^{3p}\cdot k+\epsilon ^{4p}\cdot k=0$ implies the relation
\begin{equation*}
k+(-1)^{2(p(5p+1))}\epsilon ^{p}\cdot k+(-1)^{2(2p(5p+1))}\epsilon
^{2p}\cdot k+(-1)^{2(3p(5p+1))}\epsilon ^{3p}\cdot
k+(-1)^{2(4p(5p+1))}\epsilon ^{4p}\cdot k=0
\end{equation*}%
or equivalently%
\begin{equation*}
k+(-1)^{p(5p+1)}\epsilon ^{p}\cdot _{\text{m}}k+(-1)^{2p(5p+1)}\epsilon
^{2p}\cdot _{\text{m}}k+(-1)^{3p(5p+1)}\epsilon ^{3p}\cdot _{\text{m}%
}k+(-1)^{4p(5p+1)}\epsilon ^{4p}\cdot _{\text{m}}k=0\text{.}
\end{equation*}%
Passing to coinvariants we have
\begin{equation*}
k+(-1)^{p(5p+1)}k+k+(-1)^{3p(5p+1)}k+k\sim 0~\Longleftrightarrow ~5k\sim 0.
\end{equation*}

\noindent \textbf{(D)} Let $l,k,h\in H_{2}(\cup \widehat{\mathcal{A}}(\alpha
),\mathbb{Z})$ be homology classes induced by subspaces $L$, $L\cap \epsilon
L$ and $\bigcap_{i=0}^{5}\epsilon ^{i}L$. Then $l$ generates the first copy
of $\mathbb{Z[Z}_{6}\mathbb{]}$ and $k$ the second copy of $\mathbb{Z[Z}_{6}%
\mathbb{]}$ in the $\mathbb{Z}_{6}$-module $H_{2}(\cup \widehat{\mathcal{A}}%
(\alpha ),\mathbb{Z})$. There are no relations on $l$ and $k$, so the
classes of $l$ and $k$ generate the first and the second copy of $\mathbb{Z}$
in coinvariants. The equality $\bigcap_{i=0}^{5}\epsilon ^{i}L=\epsilon
\cdot (\bigcap_{i=0}^{5}\epsilon ^{i}L)$ implies the homology equality $%
h=\epsilon \cdot h$. Thus we have the relation%
\begin{equation*}
h\sim \epsilon \cdot _{\text{m}}h=(-1)^{7}\epsilon \cdot h=-h
\end{equation*}%
which produces the $\mathbb{Z}_{2}$ in coinvariants.

\noindent \textbf{(E) }Let $l,k\in H_{n-4}(\cup \widehat{\mathcal{A}}(\alpha
),\mathbb{Z})$ be homology classes induced by subspaces $L$, $L\cap \epsilon
^{q}L$. There are no non-trivial relations on $l$ and $k$, so the
coinvariants are $\mathbb{Z\oplus Z}$.

\noindent \textbf{(F)} Let $l,k\in H_{n-4}(\cup \widehat{\mathcal{A}}(\alpha
),\mathbb{Z})$ be homology classes induced by subspaces $L$, $L\cap \epsilon
^{p+q}L$. Again, there are no relations on $l$ and $k$ and the result
follows directly.

\noindent \textbf{(G)} Let $l,k,h\in H_{n-4}(\cup \widehat{\mathcal{A}}%
(\alpha ),\mathbb{Z})$ be homology classes induced by subspaces $L$, $L\cap
\epsilon ^{p+q}L$ and $L\cap \epsilon ^{p+q}L\cap \epsilon ^{2p+2q}L$. There
is only one relevant equality $L\cap \epsilon ^{p+q}L\cap \epsilon
^{2p+2q}L=\epsilon ^{p+q}\cdot (L\cap \epsilon ^{p+q}L\cap \epsilon
^{2p+2q}L)$ which produces the equality $h=(-1)^{(p+q)(n+1)}\epsilon
^{p+q}\cdot h$ in homology (because $\det (\epsilon ^{p+q}|_{L\cap \epsilon
^{p+q}L\cap \epsilon ^{2p+2q}L})=(-1)^{(p+q)(n+1)}$). These relations become
trivial in coinvariants
\begin{equation*}
h\sim \epsilon ^{p+q}\cdot _{\text{m}}h=(-1)^{(p+q)(n+1)}\epsilon
^{p+q}\cdot h=h.
\end{equation*}

\noindent \textbf{(H) }Let $l\in H_{n-4}(\cup \widehat{\mathcal{A}}(\alpha ),%
\mathbb{Z})$ be the homology class induced by subspaces $L$. Then $%
l,\epsilon \cdot l,...,\epsilon ^{n-1}\cdot l$ are generators of $%
H_{n-4}(\cup \widehat{\mathcal{A}}(\alpha ),\mathbb{Z})$ and there are no
non trivial relations on these generators.
\end{proof}

\subsubsection{Case 2: $Q_{4n}$-coinvariants of $H_{2}(M;\mathbb{Z})$}

\begin{theorem}
\label{Th.KoinvarijanteQ4n}Let $M=W_{n}-\cup \mathcal{A}(\alpha )$ be the
complement of the minimal $D_{2n}$-arrangement $\mathcal{A}(\alpha )$
containing the subspace $L(\alpha )$ defined by (\ref{L(alfa)}) (or $Q_{4n}$%
-arrangement, where $j^{2}$ acts trivially). Let $\alpha
=(a_{1},a_{2},a_{3},a_{4})\in \mathbb{N}^{4}$ and $a_{1}+a_{2}+a_{3}+a_{4}=n$%
. Then

\begin{equation*}
\begin{tabular}{|l|l|l|}
\hline
& $\alpha $ & $H_{2}(M;\mathbb{Z})_{Q_{4n}}$ \\ \hline
(A) & $(p,p,p,p)$ & $\mathbb{Z}_{2}$ \\ \hline
(B) & $(p,q,p,q)$, $p\neq q$ & $\mathbb{Z}_{2}$ \\ \hline
(C) & $(p,p,p,2p)$ & $\mathbb{Z}_{2}\mathbb{\oplus Z}_{5}$ \\ \hline
(D) & $(1,1,1,3)$ & $\mathbb{Z}_{2}\mathbb{\oplus Z\oplus Z}_{2}$ \\ \hline
(E) & $(q,q,q,p)$, $p\notin \{q,2q\}$ & $\mathbb{Z}_{2}\mathbb{\oplus Z}$ \\
\hline
(F) & $(p,q,p,p+q)$, $p\neq q$ & $\mathbb{Z}_{2}\mathbb{\oplus Z}$ \\ \hline
(G) & $(p,p,2p,2p)$ & $\mathbb{Z\oplus Z\oplus Z}_{2}$ \\ \hline
(H) & $(p,q,p+q,p+q)$, $p\neq q$ & $?$ \\ \hline
(I) & $(p,p,q,q)$, $p\neq q$, $p\neq 2q$, $q\neq 2p$ & $\mathbb{Z}$ \\ \hline
(J) & $(p,q,r,q)$, $r\notin \{q,p,p+q\}$ & $\mathbb{Z}_{2}$ \\ \hline
(K) & $(p,q,r,p+q)$, $r\notin \{p,p+q\}$ & $?$ \\ \hline
(L) & $(p,p,q,q+r)$, $p\neq q$, $r\neq 0$ & $?$ \\ \hline
(M) & In all the remaining cases & $\mathbb{Z}$ \\ \hline
\end{tabular}%
\end{equation*}
\end{theorem}

\noindent The proof goes along the lines of the proof of theorem \ref%
{Th.KoinvarijanteZn}. Again, to compute the group  $H_{2}(M;\mathbb{Z})_{Q_{n}}$
of $Q_{4n}$-coinvariants, we use the $Q_{n}$-module $H_{n-4}(\cup \widehat{%
\mathcal{A}}(\alpha ),\mathbb{Z})$, but with the modified action. Precisely,
for $l\in H_{n-4}(\cup \widehat{\mathcal{A}}(\alpha ),\mathbb{Z})$ and $g\in
Q_{n}$,
\begin{equation*}
g\cdot _{\text{m}}l=\det (g)~g\cdot l\text{.}
\end{equation*}%
We restrict ourselves to pointing to those relations on generators of $%
H_{n-4}(\cup \widehat{\mathcal{A}}(\alpha ),\mathbb{Z})$ which produce
non-trivial identities in coinvariants. Again, let $\sim $ be the
equivalence relation on $H_{n-4}(\cup \widehat{\mathcal{A}}(\alpha ),\mathbb{%
Z})$ with the zero equivalence class generated by elements of the form $%
g\cdot _{\text{m}}x-x$, $g\in Q_{4n}$, $x\in H_{n-4}(\cup \widehat{\mathcal{A%
}}(\alpha ),\mathbb{Z})$. Thus $H_{2}(M;\mathbb{Z})_{Q_{4n}}\cong
H_{n-4}(\cup \widehat{\mathcal{A}}(\alpha ),\mathbb{Z})/\sim $.

\begin{proof}
\textbf{(A) }Since $H_{2}(M;\mathbb{Z})_{Q_{4n}}$ is a quotient of $H_{2}(M;%
\mathbb{Z})_{\mathbb{Z}_{n}}$, the result follows.

\noindent \textbf{(B) }For $\alpha =(p,q,p,q)$, let $l\in H_{n-4}(\cup
\widehat{\mathcal{A}}(\alpha ),\mathbb{Z})$ be determined by the subspace $%
L. $There are two equalities: $L=\epsilon ^{p+q}\cdot L$ and $L=\epsilon
^{-q}j\cdot L$, which in homology become equalities $l=(-1)^{p+q}\epsilon
^{p+q}\cdot l$ and $l=(-1)^{p+2q+1}\epsilon ^{-q}j\cdot l$. These equalities
imply that
\begin{equation*}
l\sim \epsilon ^{p+q}\cdot _{\text{m}}l=\det (\epsilon ^{p+q})~\epsilon
^{p+q}\cdot l=(-1)^{(p+q)(n+1)}~\epsilon ^{p+q}\cdot
l=(-1)^{(p+q)(2p+2q+1)}(-1)^{p+q}l=l,
\end{equation*}%
and%
\begin{equation*}
l\sim \epsilon ^{-q}j\cdot _{\text{m}}l=\det (\epsilon ^{-q}j)~\epsilon
^{-q}j\cdot l=(-1)^{(n+1)q+\binom{{n}}{{2}}}~\epsilon ^{-q}j\cdot
l=(-1)^{p+2q}(-1)^{p+2q+1}l=-l\text{.}
\end{equation*}

\noindent \textbf{(C)} Let $l,k\in H_{n-4}(\cup \widehat{\mathcal{A}}(\alpha
),\mathbb{Z})$ be homology classes induced by subspaces $L$ and $L\cap
\epsilon ^{p}L$. In contrast to the proof of theorem \ref{Th.KoinvarijanteZn}
(C), we have a relation on $l$, which is a consequence of the equality $%
L=\epsilon ^{-2p}j\cdot L$. Thus in homology $l=(-1)^{\binom{{5p}}{{2}}%
+1}\epsilon ^{-2p}j\cdot l$ and consequently in coinvariants%
\begin{equation*}
l\sim \epsilon ^{-2p}j\cdot _{\text{m}}l=(-1)^{2p(5p+1)+\binom{{5p}}{{2}}%
}\epsilon ^{-2p}j\cdot l=(-1)^{\binom{{5p}}{{2}}}\epsilon ^{-2p}j\cdot l=-l.
\end{equation*}%
There is also one more relation on $k$ which is the consequence of the
equality $L\cap \epsilon ^{p}L=\epsilon ^{-p}j\cdot (L\cap \epsilon ^{p}L)$.
Passing to homology we get the equality $k=(-1)^{p(n+1)+\binom{n}{{2}}%
}\epsilon ^{-p}j\cdot k$. In coinvariants we get the trivial relation%
\begin{equation*}
k\sim \epsilon ^{-p}j\cdot _{\text{m}}k=(-1)^{p(n+1)+\binom{n}{{2}}}\epsilon
^{-p}j\cdot k=k
\end{equation*}%
and the second summand is just $\mathbb{Z}_{5}$.

\noindent \textbf{(D)} Let $l,k,h\in H_{2}(\cup \widehat{\mathcal{A}}(\alpha
),\mathbb{Z})$ be induced by subspaces $L$, $L\cap \epsilon L$ and $%
\bigcap_{i=0}^{5}\epsilon ^{i}L$. The equalities $L=\epsilon ^{-3}j\cdot L$
and $L\cap \epsilon L=\epsilon ^{-2}j\cdot (L\cap \epsilon L)$ imply the
equalities $l=-\epsilon ^{-3}j\cdot l$ and $k=-\epsilon ^{-2}j\cdot k$
(since $\det (\epsilon ^{-3}j|_{L})=-1$ and $\det (\epsilon ^{-2}j|_{L\cap
\epsilon L})=-1$). We do not consider equality $\bigcap_{i=0}^{5}\epsilon
^{i}L=j(\bigcap_{i=0}^{5}\epsilon ^{i}L)$, because $\bigcap_{i=0}^{5}%
\epsilon ^{i}L=\epsilon (\bigcap_{i=0}^{5}\epsilon ^{i}L)$ already produced $%
\mathbb{Z}_{2}$-torsion. Thus in coinvariants%
\begin{equation*}
l\sim \epsilon ^{-3}j\cdot _{\text{m}}l=\epsilon ^{-3}j\cdot l=-l\text{ \
and }k\sim \epsilon ^{-2}j\cdot _{\text{m}}k=-\epsilon ^{-2}j\cdot k=k\text{.%
}
\end{equation*}

\noindent \textbf{(E) }Let $l,k\in H_{n-4}(\cup \widehat{\mathcal{A}}(\alpha
),\mathbb{Z})$ be homology classes induced by subspaces $L$, $L\cap \epsilon
^{q}L$. Now we have two equalities $L=\epsilon ^{-p}j\cdot L$ and $L\cap
\epsilon ^{q}L=\epsilon ^{-(2p+2q)}j\cdot (L\cap \epsilon ^{q}L)$ in
contrast to the similar case in the proof of the theorem \ref%
{Th.KoinvarijanteZn}, (E).These equalities produce equalities $%
l=(-1)^{(n+1)3q+\binom{n}{{2}}+1}\epsilon ^{-p}j\cdot l$ and $%
k=(-1)^{(n+1)(q-p)+\binom{n}{{2}}}\epsilon ^{q-p}j\cdot k$ in homology.
Passing to coinvariants we have%
\begin{equation*}
l\sim \epsilon ^{-p}j\cdot _{\text{m}}l=(-1)^{(n+1)3q+\binom{n}{{2}}%
}\epsilon ^{-p}j\cdot l=-l\text{ and }k\sim \epsilon ^{q-p}j\cdot _{\text{m}%
}k=(-1)^{(n+1)(q-p)+\binom{n}{{2}}}\epsilon ^{q-p}j\cdot k=k\text{.}
\end{equation*}

\noindent \textbf{(F) }Let $l,k\in H_{n-4}(\cup \widehat{\mathcal{A}}(\alpha
),\mathbb{Z})$ be homology classes induced by subspaces $L$, $L\cap \epsilon
^{p+q}L$. The equalities $L=\epsilon ^{-(p+q)}j\cdot L$ and $L\cap \epsilon
^{p+q}L=j\cdot (L\cap \epsilon ^{p+q}L)$ imply $l=(-1)^{(n+1)(2p+q)+\binom{n%
}{{2}}+1}\epsilon ^{-(p+q)}j\cdot l$ and $k=(-1)^{\binom{n}{{2}}}j\cdot k$
in homology. Therefore, in coinvariants%
\begin{equation*}
l\sim \epsilon ^{-(p+q)}j\cdot _{\text{m}}l=(-1)^{(n+1)(2p+q)+\binom{n}{{2}}%
}\epsilon ^{-(p+q)}j\cdot l=-l\text{ \ and \ }k\sim j\cdot _{\text{m}%
}k=(-1)^{\binom{n}{{2}}}j\cdot k=k\text{.}
\end{equation*}

\noindent \textbf{(G)} Let $l,k,h\in H_{n-4}(\cup \widehat{\mathcal{A}}%
(\alpha ),\mathbb{Z})$ be homology classes induced by subspaces $L$, $L\cap
\epsilon ^{4p}L$ and $L\cap \epsilon ^{2p}L\cap \epsilon ^{4p}L$. The
relevant equalities are%
\begin{equation*}
\epsilon ^{2p}j\cdot L=L\text{, }j\cdot (L\cap \epsilon ^{4p}L)=L\cap
\epsilon ^{4p}L\text{, }j\cdot (L\cap \epsilon ^{2p}L\cap \epsilon
^{4p}L)=L\cap \epsilon ^{2p}L\cap \epsilon ^{4p}L
\end{equation*}%
which in homology imply%
\begin{equation*}
l=(-1)^{\binom{n}{{2}}}\epsilon ^{2p}j\cdot l\text{, }k=(-1)^{\binom{n}{{2}}%
}j\cdot k\text{, }h=(-1)^{\binom{n}{{2}}+1}j\cdot h\text{.}
\end{equation*}%
So in coinvariants we have the following relations%
\begin{equation*}
l\sim \epsilon ^{2p}j\cdot _{\text{m}}l=(-1)^{\binom{n}{{2}}}\epsilon
^{2p}j\cdot l=l\text{, }k\sim j\cdot _{\text{m}}k=k=(-1)^{\binom{n}{{2}}%
}j\cdot k=k\text{, }h\sim j\cdot _{\text{m}}h=(-1)^{\binom{n}{{2}}}j\cdot
h=-h\text{.}
\end{equation*}

\noindent \textbf{(I)} Let $l\in H_{n-4}(\cup \widehat{\mathcal{A}}(\alpha ),%
\mathbb{Z})$ be the homology class induced by the subspace $L$. There is
only one equality, $L=\epsilon ^{-2q}j\cdot L$ which in homology reads as $%
l=(-1)^{\binom{n}{{2}}}\epsilon ^{-2q}j\cdot l$ and in coinvariants%
\begin{equation*}
l\sim \epsilon ^{-2q}j\cdot _{\text{m}}l=(-1)^{\binom{n}{{2}}}\epsilon
^{-2q}j\cdot l=l
\end{equation*}

\noindent \textbf{(J) }Let $l\in H_{n-4}(\cup \widehat{\mathcal{A}}(\alpha ),%
\mathbb{Z})$ be induced by the subspace $L$. The equality $\epsilon
^{p}j\cdot L=L$ implies equality $l=(-1)^{p(n+1)+\binom{n}{{2}}+1}\epsilon
^{p}j\cdot l$ in homology. As before, in coinvariants this leads to the
equality
\begin{equation*}
l\sim \epsilon ^{p}j\cdot _{\text{m}}l=(-1)^{p(n+1)+\binom{n}{{2}}}\epsilon
^{p}j\cdot l=-l\text{.}
\end{equation*}

\noindent \textbf{(M)} In all the remaining cases there are no relations
involving the element $l\in H_{n-4}(\cup \widehat{\mathcal{A}}(\alpha ),%
\mathbb{Z})$ induced by subspace $L$. Since the homology group $H_{n-4}(\cup
\widehat{\mathcal{A}}(\alpha ),\mathbb{Z})$ is freely generated by the orbit
of $l$, the group $H_{n-4}(\cup \widehat{\mathcal{A}}(\alpha ),\mathbb{Z}%
)_{Q_{4n}}$ of coinvariants is $\mathbb{Z}$.
\end{proof}

\subsection{The main theorem}

\noindent This is the central section of the paper. We gather together all
the information collected in previous sections and use it to obtain a
reasonably complete answer to the Problem \ref{prob:dihedral}.

\begin{theorem}
There does not exist a $\mathbb{Z}_{n}\text{-map
}F:S^{3}\rightarrow M(\alpha )$ if and only if
\label{Th.Main1}%
\begin{equation*}
(\exists \text{ }p\in \mathbb{N})\text{ }%
\alpha =(p,p,p,p)\text{ or }\alpha =(p,p,p,2p)\text{.}
\end{equation*}
\end{theorem}

\begin{proof}
(A) We first deal with the existence of a $\mathbb{Z}_{n}$-map. Assume that $%
\alpha $ is not of the form $(p,p,p,p)$ or $(p,p,p,2p)$. Then the cohomology
class of the obstruction cocycle $[c_{\mathbb{Z}_{n}}^{\prime }(h)(d_{3})]$
is divisible by $2$ (Theorem \ref{Th.ObCycleZn}). On the other hand this
cohomology class is a torsion element of $H_{2}(M;\mathbb{Z})_{\mathbb{Z}%
_{n}}$ (Proposition \ref{Th.ObCoCycleTorzioni}). Since Theorem \ref%
{Th.KoinvarijanteZn} implies that the only torsion summands which can appear
are copies of $\mathbb{Z}_{2}$, we conclude that $[c_{\mathbb{Z}%
_{n}}^{\prime }(h)(d_{3})]=0$. Thus a $\mathbb{Z}_{n}$-map $%
F:S^{3}\rightarrow M(\alpha )$ must exists.

\noindent (B) Let $\alpha =(p,p,p,p)$. In order to prove that $[c_{\mathbb{Z}%
_{n}}^{\prime }(h)(d_{3})]=(-1)^{p}{[[}y_{1}]]+K\neq 0$ we again rely on the
isomorphism $\varphi :H_{2}(M(\mathcal{A}(\alpha )),\mathbb{Z})\rightarrow
\mathrm{Hom}(H_{n-4}(\cup \widehat{\mathcal{A}}(\alpha ),\mathbb{Z}),\mathbb{%
Z})$. Using the fact that this isomorphism is related to the linking number,
we have that
\begin{equation*}
{[[}y_{1}]]\rightarrow \tsum\limits_{i=0}^{p-1}\mathrm{link}(\widehat{%
\epsilon ^{i}L},\widehat{(y_{1}+L^{\bot })})\epsilon ^{i}l=l
\end{equation*}%
since $\widehat{L}$ links only with $\widehat{(y_{1}+L^{\bot })}$. Here $%
\widehat{\epsilon ^{i}L}$ and $\widehat{(y_{1}+L^{\bot })}$ are respectively
one-point compactifications of the linear subspace $\epsilon ^{i}L$ and the
affine subspace $y_{1}+L^{\bot }$ in the sphere $\widehat{W_{n}}\approx
S^{n-1}$. Now let $K$ be the additive subgroup of $H_{2}(M;\mathbb{Z})$
generated by the elements of the form $g\cdot x-x$, $g\in \mathbb{Z}_{n}$, $%
x\in H_{2}(M;\mathbb{Z})$. Then the group of coinvariants is $H_{2}(M;%
\mathbb{Z})_{\mathbb{Z}_{n}}=H_{2}(M;\mathbb{Z})/K$. Since $[c_{\mathbb{Z}%
_{n}}^{\prime }(h)(d_{3})]=(-1)^{p}{[[}y_{1}]]+K$ we have that $[c_{\mathbb{Z%
}_{n}}^{\prime }(h)(d_{3})]$ is the generator in $H_{2}(M;\mathbb{Z})_{%
\mathbb{Z}_{n}}\cong \mathbb{Z}_{2}$. We conclude that there does not exist
a $\mathbb{Z}_{n}$-map $F:S^{3}\rightarrow M(\alpha )$.

\noindent (C) Let $\alpha =(p,p,p,2p)$. Like in the preceding case we look
at the $\varphi $ image of the point classes ${[[}y_{1}]]$ and ${[[}y_{2}]]$,%
\begin{eqnarray*}
{[[}y_{1}]] &\rightarrow &\tsum\limits_{i=0}^{5p-1}\mathrm{link}(\widehat{%
\epsilon ^{i}L},\widehat{(y_{1}+L^{\bot })})\epsilon
^{i}l+\tsum\limits_{i=0}^{4p-1}\mathrm{link}(\widehat{\epsilon ^{i}(L\cap
\epsilon ^{p}L)},\widehat{(y_{1}+L^{\bot })})\epsilon ^{i}k=l+k \\
{[[}y_{2}]] &\rightarrow &\tsum\limits_{i=0}^{5p-1}\mathrm{link}(\widehat{%
\epsilon ^{i}L},\widehat{(y_{1}+L^{\bot })})\epsilon
^{i}l+\tsum\limits_{i=0}^{4p-1}\mathrm{link}(\widehat{\epsilon ^{i}(L\cap
\epsilon ^{p}L)},\widehat{(y_{1}+L^{\bot })})\epsilon ^{i}k=l\text{.}
\end{eqnarray*}%
Since in this case $[c_{\mathbb{Z}_{n}}^{\prime
}(h)(d_{3})]=2[[y_{1}]]-2[[y_{2}]]+K$ we have that $[c_{\mathbb{Z}%
_{n}}^{\prime }(h)(d_{3})]\in H_{2}(M;\mathbb{Z})_{\mathbb{Z}_{n}}\cong
\mathbb{Z\oplus Z}_{5}$ is $2\cdot ($generator of the second summand$)$ and
so must be different from zero. Therefore, there are no $\mathbb{Z}_{n}$%
-maps $F:S^{3}\rightarrow M(\alpha )$.
\end{proof}

\begin{theorem}
\label{Th.Main2}A partial answer to the problem \ref{prob:dihedral} for the
whole group $Q_{4n}$ is the following.

(A) If $\alpha =(p,p,p,p)$ or $\alpha =(p,p,p,2p)$, then there are no $%
Q_{4n} $-maps $F:S^{3}\rightarrow M(\alpha )$.

(B) If $\alpha $ is not of the form $(p,q,p+q,p+q)$, $p\neq q$, or $%
(p,q,r,p+q)$, $r\notin \{p,p+q\}$, or $(p,p,q,q+r)$, then there exists a $%
Q_{4n}$-map $F:S^{3}\rightarrow M(\alpha )$.

(C) If $\alpha =$ $(p,q,p+q,p+q)$, $p\neq q$, or $\alpha =(p,q,r,p+q)$, $%
r\notin \{p,p+q\}$, or $\alpha =(p,p,q,q+r)$, the methods of the paper are
inconclusive and a further analysis is needed.
\end{theorem}

\begin{proof}
The proof goes along the lines of the proof of the preceding theorem. The
case (A) is the direct consequence of the preceding theorem, but it can be
also derived by a direct calculation interpreting the $Q_{4n}$ obstruction
cocycle in the appropriate group of coinvariants. The case (B) follows from
the following facts:

(1) the cohomology class of the obstruction cocycle $[c_{Q_{4n}}^{\prime
}(h)(d_{3})]$ is divisible by $4$ (Theorem \ref{Th.ObCycleQ4n}),

(2) the cohomology class of the obstruction cocycle $[c_{Q_{4n}}^{\prime
}(h)(d_{3})]$ is a torsion element (Proposition \ref{Th.ObCoCycleTorzioni}),
and

(3) the only torsion group which appears as a summand in the coinvariant
group is $\mathbb{Z}_{2}$ (Theorem \ref{Th.KoinvarijanteQ4n}).
\end{proof}

\section{{Appendix}}

\subsection{The geometry of the group $Q_{4n}$}

\label{sec:geometry}

\noindent It is of utmost importance for computations to describe an
economical and geometrically transparent $G$-invariant, CW-structure on a
given $G$-manifold. In the case of the generalized quaternion group $Q_{4n}$
acting on the sphere $S^{3}$ of all unit quaternions, such a structure is
described by the following construction.

\begin{figure}[htb]
\centering
\includegraphics[scale=0.45]{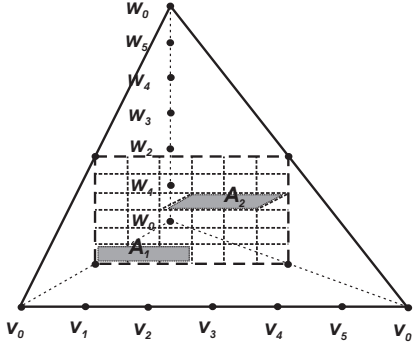}
\caption{The sphere $S^3 = P^{(1)}_{2n}\ast P^{(2)}_{2n}$ in the
case $n=3$} \label{fig:taraba2}
\end{figure}

\noindent We have already met in Section~\ref{sec:obs-theo} the join
decomposition $S^{3}=P_{2n}^{(1)}\ast P_{2n}^{(2)}$ of the $3$-sphere. The
case $n=3$ is depicted on the Figure \ref{fig:taraba2} where $%
S^{3}=S^{1}\ast S^{1}=[0,1]\ast \lbrack 0,1]/\approx $ and \textquotedblleft
$\approx $\textquotedblright\ is the equivalence relation on the $3$-simplex
$\sigma ^{3}:=[0,1]\ast \lbrack 0,1]$ arising from the identification $%
S^{1}\cong \lbrack 0,1]/0\sim 1$. The (big) rectangle represents the torus $%
T^{2}:=\{\frac{1}{2}x+\frac{1}{2}y\mid x,y\in S^{1}\}$ subdivided into small
rectangles which are in $1$--$1$ correspondence with simplices $%
[v_{i},v_{i+1};w_{j},w_{j+1}]$. Since each point $a=\frac{1}{2}x+\frac{1}{2}%
y $ of this rectangle uniquely determines the line segment $[x,y]\subset
S^{3}$, each subset $A\subset T^{2}$ determines uniquely a set $\hat{A}%
\subset S^{3}$ such that $\hat{A}\cap T^{2}=A$. Let $A_{1}:=\cup
_{j=0}^{n-1}[v_{i},v_{i+1}]\times \lbrack w_{0},w_{1}]$ be a rectangle in $%
T^{2}$ and $A_{2}$ a closely related \textquotedblleft
skew\textquotedblright\ rectangle, see the Figure
\ref{fig:taraba2}. Then
both $\hat{A}_{1}$ and $\hat{A}_{2}$ are fundamental domains for the $Q_{4n}$%
-action on $S^{3}$, both inducing $Q_{4n}$-invariant $CW$-structures on $%
S^{3}$. The second of these two structures is simpler having only one $%
Q_{4n} $-cell $a$ in dimension $0$, two $Q_{4n}$-cells $b$ and $b^{\prime }$
in dimension $1$, two $Q_{4n}$-cells $c$ and $c^{\prime }$ in dimension $2$,
and one $Q_{4n}$-cell $e$ in dimension $3$. The reader can easily check that
the associated cellular chain complex $\{D_{i}(S^{3},\mathbb{Z})\}_{i=0}^{3}$
of $Q_{4n}$-modules coincides with the well known minimal resolution of $%
\mathbb{Z}$ by free $Q_{4n}$-modules described on p.~253 of the classical
monograph \cite{CaEi}. More precisely
\begin{equation*}
0\rightarrow Z(Q_{4n})e\overset{\partial }{\rightarrow }Z(Q_{{4n}})c\oplus
Z(Q_{{4n}})c^{\prime }\overset{\partial }{\rightarrow }Z(Q_{{4n}})b\oplus
Z(Q_{{4n}})b^{\prime }\overset{\partial }{\rightarrow }Z(Q_{4n})e\rightarrow
0
\end{equation*}%
where
\begin{equation*}
\begin{array}{lll}
\partial e=(\epsilon -1)c-(\epsilon j-1)c^{\prime }; & \partial
c=Nb-(j+1)b^{\prime }; & \partial c^{\prime }=(\epsilon j+1)b+(\epsilon
-1)b^{\prime } \\
& \partial b=(\epsilon -1)a & \partial b^{\prime }=(j-1)a%
\end{array}%
\end{equation*}%
and $N=1+\epsilon +\ldots +\epsilon ^{n-1}\in \mathbb{Z}(Q_{4n})$. On
applying the functor $\mathrm{Hom}_{Q_{4n}}(\cdot ,N))$, where $N:=H_{2}(M;%
\mathbb{Z})$ is a $Q_{4n}$-module, one obtains the sequence,
\begin{equation*}
0\longleftarrow H_{2}(M;\mathbb{Z})\overset{\Gamma }{\longleftarrow }H_{2}(M;%
\mathbb{Z})\oplus H_{2}(M;\mathbb{Z})\leftarrow H_{2}(M;\mathbb{Z})\oplus
H_{2}(M;\mathbb{Z})\longleftarrow H_{2}(M;\mathbb{Z})\longleftarrow 0
\end{equation*}%
where $\Gamma (p,q)=(\epsilon -1)p-(\epsilon j-1)q$ for $p,q\in H_{2}(M;%
\mathbb{Z})$. Now it is not hard to prove
\begin{equation*}
H_{Q_{4n}}^{3}(S^{3},H_{2}(M;\mathbb{Z}))=H_{2}(M;\mathbb{Z})/\mathrm{Im}%
\Gamma \cong H_{2}(M;\mathbb{Z})_{Q_{4n}}
\end{equation*}%
where $H_{2}(M;\mathbb{Z})_{Q_{4n}}$ is a group of coinvariants of the $%
Q_{4n}$-module $H_{2}(M;\mathbb{Z})$, \cite{Brown}. Alternatively, the last
result can be seen as a consequence of the equivariant Poincar\'{e} duality,
\cite{Wall}.

\noindent The complex associated with the fundamental domain $A_{1}$, Figure~%
\ref{fig:taraba2}, is not as economical as the complex based on
the fundamental domain $A_{2}$. Note however that for our purposes
this is not a
serious disadvantage. Indeed, all we need to know is the intersection of $%
e_{1}=\mathrm{int}(A_{1})$ with the singular set $h^{-1}(\cup \mathcal{A}%
(\alpha ))$, for a generic $Q_{4n}$-equivariant map $h$, and
interpret the answer as an element of the group
$H_{2}(M;\mathbb{Z})_{Q_{4n}}$, cf. Sections \ref{sec:genpos} and
\ref{sec:obscoc}. On the other hand, the choice of $A_{1}$ has an
advantage that its top cell $e_{1}=A_{1}$ is a union of simplices
from the original simplicial structure on $S^{3}$, see
Figure~\ref{fig:taraba2}, which considerably simplifies the
computations performed in Section~\ref{sec:obscoc}.

\subsection{Change of the group}

\noindent Let $X$ be a left $G$-space and $H\vartriangleleft G$ be a normal
subgroup. Then $X/H$ is a $G/H$-space with the action given by $gH(Hx)=H(gx)$%
. A $G/H$-space $Z$ can be also seen as a $G$-space via the quotient
homomorphism $\pi :G\rightarrow G/H$ where for $g\in G$ and $z\in Z,\,g\cdot
z:=\pi (g)z$.

\begin{proposition}
\label{PropositionScalar} Suppose that $X$ and $Z$ are $G$-spaces and $%
H\vartriangleleft G$ is a normal subgroup of $G$ which acts trivially on $Z$%
. Then there exists a $G$-map $\alpha :X\rightarrow Z$ if and only if there
exists a $G/H$-map $\beta :X/H\rightarrow Z$, where $X/G$ and $Z/H=Z$ are
interpreted as $G/H$-spaces.
\end{proposition}

\begin{proof}
$\Leftarrow :$ Observe that the quotient map $p:X\rightarrow X/H$ is a $G$%
-map and that the $G/H$-space $X/H$ is a $G$-space via the homomorphism $\pi
:G\rightarrow G/H$. If $\beta :X/H\rightarrow Z$ is a $G/H$-map, it is also
a $G$-map and the composition $\alpha :=p\circ \beta :X\rightarrow Z$ is a $%
G $-map.

\noindent $\Rightarrow :$ Let $\alpha :X\rightarrow Z$ be a $G$-map. Then,
since $H$-acts trivially on $Z$, there is a factorization $\alpha =\beta
\circ p$ for some $\beta :X/H\rightarrow Z$. We check that $\beta $ is a $%
G/H $-map,
\begin{equation*}
\beta (gH\cdot Hx)=\alpha (gx)=g\alpha (x)=g\beta (Hx)=(gH)\beta (Hx).
\end{equation*}
\end{proof}

\subsection{Modifying free $G$-action on $S^{3}$}

\label{sec:modifying}

\noindent

\label{PropositionChangeofAction}(A) Suppose that $\gamma _{i}:G\times
S^{3}\rightarrow S^{3},\gamma _{i}(g,x)=g\cdot _{i}x,\,i=1,2$ are two
actions of a finite group $G$ on the $3$-sphere $S^{3}$, and assume that the
action $\gamma _{1}$ is free. Then there exists a $G$-equivariant map $%
f:S^{3}\rightarrow S^{3}$ between these two actions, i.e.
\begin{equation*}
(\forall g\in G)\,(\forall x\in S^{3})\,\,f(g\cdot _{1}x)=g\cdot _{2}f(x).
\end{equation*}%
(B) If in addition the action $\gamma _{2}$ is free, then for any $G$-space $%
Z$ there exists a $\gamma _{1}$-map $f:S^{3}\rightarrow Z$ if and only if
such a map exists for the $\gamma _{2}$ action.

\begin{proof}
(A) The proof is routine and relies on the fact that $S^{3}$ is a $2$%
-connected, $\gamma _{1}$-free, $CW$-complex so there are no obstructions to
extend equivariantly a map defined on the $0$-skeleton of $S^{3}$. The
statement (B) is a direct consequence of (A).
\end{proof}

\end{document}